\documentclass[11pt]{article}
\usepackage{a4}
\usepackage{graphicx}
\usepackage{parskip}

\usepackage{amsfonts}
\usepackage{natbib}
\usepackage{algorithm}
\usepackage{algorithmic}

\usepackage{latexsym}

\newcommand{\R}{\mathbb{R}}

\newcommand{\PP} {{  \rm I\hskip-0.22em P}}

\newcommand{\EE} {{\rm I\hskip-0.48em E}}

\usepackage{natbib}

\usepackage{amsmath} 
\usepackage{amsthm}
\usepackage{amssymb} 
\usepackage{mathrsfs} 
\usepackage{amsfonts} 
\usepackage{mathtools}
\usepackage{algorithm} 
\usepackage{algorithmic} 
\usepackage{refcount}
\usepackage{caption} 
\usepackage{subcaption} 
\usepackage{graphicx}
\usepackage{float}
\usepackage[bottom]{footmisc}

\newtheorem{theorem}{Theorem}[section]

\newtheorem{lemma}{Lemma}[section]
\newtheorem{corollary}{Corollary}[section]
\newtheorem{remark}{Remark}[section]
\newtheorem{example}{Example}[section]
\newtheorem{condition}{Condition}[section]

\begin{document}

\centerline{\bf \Large The additive model}

\centerline{\bf \Large with different smoothness for the components}

\vskip .1in
\centerline{Sara van de Geer and Alan Muro}

\centerline{Seminar for Statistics, ETH Z\"urich}

\vskip .1in
\centerline {May 26, 2014} 

{\bf Abstract} We consider an additive regression model consisting of two components
$f^0$ and $g^0$, where the first component $f^0$ is in some sense ``smoother" than the second $g^0$.
Smoothness is here described in terms of a semi-norm on the class of regression functions.
We use a penalized least squares estimator $(\hat f , \hat g )$ of $(f^0 , g^0 )$ and
show that the rate of convergence for $\hat f $ is faster than the rate of convergence
for $\hat g$. In fact, both rates are generally as fast as in the case where one
of the two components is known. The theory is illustrated by a
simulation study. Our proofs rely on recent results from empirical
process theory. 



\section{Introduction}\label{introduction.section}
Additive modelling  has a long history (\cite{stone1985additive},
 \cite{hastie1990generalized}) and is very useful for
 dealing with the curse of dimensionality. Important estimation methods
 for such models are for example spline smoothing
 (\cite{wahba1990spline}) or iterative back fitting (\cite{mammen1999existence}). 
 Our contribution in this paper is to show that standard spline smoothing or more generally
 penalized least squares can estimate ``smoother" components
 at a faster rate than ``rough" components. In fact, we show an oracle rate for the
 smoother component, which is as fast as in the case where the rough component is
 known. Similarly (but perhaps less surprisingly) the rougher component can be
 estimated as fast as in the case where the smooth component is known. 
 These results are in the same spirit as  results for semi-parametric models
 (\cite{bicketal98}) saying that the parametric
 part (the parameter of interest) is estimated with parametric rate despite the presence of an infinite-dimensional
 nuisance parameter. We make use of recent empirical process theory to deal with an infinite-dimensional parameter of interest.
  
 For simplicity we consider the additive  model with two components (extensions to more 
 components can be derived essentially along the same lines).
Let $(X_i, Z_i)_{i=1}^n$ be i.i.d.\  input variables and $\{ Y_i \}_{i=1}^n $ be i.i.d.\ real-valued output variables.
The model is
$$Y_i = f^0 (X_i) + g^0 (Z_i) + \epsilon_i, \ i=1 , \ldots, n , $$
where $f^0 \in {\cal F}$, $g^0 \in {\cal G}$ with
${\cal F}$ and ${\cal G}$ linear function spaces. Moreover, $\epsilon:= (\epsilon_1 , \ldots , \epsilon_n)^T$ is a vector of
i.i.d.\ centered noise variables, independent of $\{ (X_i,Z_i)\}_{i=1}^n$. 
For a vector $v \in \R^n $ we write $\| v \|_n^2 := v^T v / n $. 
We study the estimator
$$(\hat f , \hat g ) := \arg \min_{f,g}\biggl \{  \| Y - f -g \|_n^2 + \lambda^2 I^2 (f) + \mu^2 J^q (g)\biggr \}  , $$
where $I$ is a semi-norm on ${\cal F}$, $J$ is a semi-norm on ${\cal G} $ and
$\lambda$ and $\mu$ are tuning parameters.  Moreover, $1 \le q \le 2 $
is some fixed constant.  
We consider the case where the ``smoothness" induced by $I$ is larger than the ``smoothness" induced by $J$.
For example, when both $X $ and $Z$ are bounded real-valued random
variables, one may think of $I$ as being some Sobolev norm, $J$ being the total
variation norm and $q=1$.  Note that we restrict ourselves to a squared norm in the penalty
for the smoother part. A generalization here is straightforward but technical. 
Also a generalization to values of $q>2$ is not difficult but is omitted to avoid complicated expressions.

We show that with an appropriate choice of the regularization parameters $\lambda$ and $\mu$ the 
rate of convergence for the smoother function $f^0$ is faster than the rate for the less smooth function $g^0$.
For each component we obtain the rate of convergence corresponding to the situation were
the other component  is known. This result is established 
assuming an incoherence condition
between $X_1$ and $Z_1$ (see Condition \ref{incoherence.condition}).

The results in this paper are related to \cite{Wahl2014}. The latter studies an additive two-component model
and applies restricted least squares instead of the penalized least squares used here.
Another important paper on the topic is \cite{efromovich2013nonparametric} where adaptive rates
are derived using a method including blockwise shrinkage.
Related is also the paper
\cite{Mueller2013} where a partial linear model is studied with the linear part being high-dimensional. 
The method used there is penalized least squares with $\ell_1$-penalty on the linear part.

\subsection{Organization of the paper} In the next section we outline the conditions used.
Main condition is an entropy condition (Condition \ref{entropy.condition}) which describes the assumed roughness
of the functions $f^0$ and $g^0$. Section \ref{main.section} contains the main theoretical result in Theorem
\ref{main.theorem}. Section \ref{simulation.section} presents a simulation study. All proofs are
in Section \ref{proofs.section}.

\section{Conditions}\label{conditions.section}
Let $P$ be the distribution of $(X,Z)$ and $\| \cdot \|$ be the $L_2 ( P)$-norm. 
For arbitrary positive constants $R$ and $M$ we let
${\cal F} (R,M):= \{ f \in {\cal F} : \ \| f \| \le R , \ I(f) \le M \} $ and $ {\cal G} ( R,M) := 
\{ g \in {\cal G}: \ \| g \| \le R , \ J(g) \le M \} $. 

Let $\| \cdot \|_{\infty} $ be the supremum norm. The entropy of $({\cal F}(R,M) , \| \cdot \|_{\infty})$
is denoted by ${\cal H}_{\infty} (\cdot , {\cal F}(R,M)) $.
The entropy integral ${\cal J}_{\infty}( \cdot  , {\cal F} (R,M) )$ is defined as
$${\cal J}_{\infty}(z, {\cal F}(R,M)) := z  \int_{0}^1 \sqrt {{\cal H}_{\infty}  ( uz, {\cal F}(R,M), \| \cdot \|_{\infty}  )}  du , \ z > 0 $$
which we assume to exist. 

For the class ${\cal G}$ the entropy ${\cal H}_{\infty} ( \cdot , {\cal G} (R,M) )$ and entropy integral
${\cal J}_{\infty} ( \cdot , {\cal G} (R,M))$ are defined similarly. 
We shall however use a somewhat relaxed version of entropy and entropy integral for ${\cal G}$.
Let ${\cal A}_n$ be the set of all subsets of cardinality $n$ within the support of $Z_1$ (equal points are allowed).
For $A_n \in {\cal A}_n$ and $g$ a real-valued function on this support we let
$$ \| g \|_{A_n, {\infty}} := \max_{z\in A_n } | g(z) | . $$
The entropy of the class $( {\cal G} (R,M), \| \cdot \|_{A_n}) $ endowed with
$\| \cdot \|_{A_n}$-norm is denoted by
${\cal H}_{A_n} ( \cdot , {\cal G} (R,M)  )$. The uniform entropy is
$${\cal H}_n ( \cdot , {\cal G} (R,M) ) := \sup_{A_n \in {\cal A}_n } {\cal H}_{A_n} ( \cdot , {\cal G} (R,M)) 
. $$
We furthermore define the entropy integral
\begin{equation}\label{Jinfty.definition}
{\cal J}_n ( z , {\cal G} (R,M) ) := z  \int_{0}^1 \sqrt {{\cal H}_n  ( uz, {\cal G}(R,M) )}  du  , \ z >0   
\end{equation}
assuming again it exists.
Note that ${\cal H}_n (\cdot , {\cal G} (R,M)) \le {\cal H}_{\infty} ( \cdot , {\cal G}(R,M)) $ and
consequently ${\cal J}_n ( \cdot , {\cal G}(R,M)) \le {\cal J}_{\infty} ( \cdot , {\cal G} (R,M)) $. 

We fix the ``roughness indices" $0< \alpha < \beta < 1 $ and assume the following bounds on the 
entropy integrals for ${\cal F}(R,M)$ and ${\cal G}(R,M)$. The reason for the more stringent version
of entropy (or entropy integral) for ${\cal F}(R,M)$ is apparent from Lemma \ref{entropy.lemma} where we consider
for $f \in {\cal F}(R,M)$
conditional versions of $f(X_1)$ given $Z_1$.

\begin{condition} \label{entropy.condition} For $R \le M$ and some constants $A_I\ge 1 $ and $A_J \ge 1$, it holds that
$${\cal J}_{\infty} ( z, {\cal F} (R,M)  ) \le A_I M^{\alpha} z^{1- \alpha}  , \ z > 0 , $$
and
$${\cal J}_n ( z, {\cal G} (R,M)  ) \le A_J M^{\beta} z^{1- \beta}  , \ z > 0 . $$
\end{condition}

As an illustration, suppose that $X_1\in [0,1]$ and
$I^2(f) = \int | f^{(k)} (x) |^2 dx $,  where $f^{(k)}$ denotes the $k$-th
derivative of $f$.   Then $\alpha=1/(2k)$ and the constant $A_I$ depends only
on the smallest eigenvalue of the matrix
$\EE \psi^T(X_1) \psi (X_1)  $ where $\psi(X_1) = ( 1, X_1 , \ldots , X_1^{k-1})$ (see
e.g.\ \cite{Birman:67}). Similar bounds hold for a general
class of Besov spaces, see
\cite{Birge:00}. 


We assume $\sup \{  \| f \|_{\infty} :\ f \in {\cal F} (R,M) \}$ is bounded by a constant proportional
to $M$ and similarly for ${\cal G} (R,M)$. Without loss of generality we assume the proportionality
constant to be equal to 1.

\begin{condition} \label{supnorm.condition} For some constant $B\ge 1$ and all $M > 0$ and any $R \le M/B$ 
it holds that $$  \sup_{f \in {\cal F} (R,M) }  \| f \|_{\infty} \le  M   , $$
and
$$ \sup_{g \in {\cal G} (R, M)  } \| g \|_{\infty} \le  M   . $$
\end{condition}

For a sub-Gaussian random variable $Z \in \R$ and $\Psi(z) := \exp [ |z|^2 ] -1$,  we define the
Orlicz norm
$$\| Z \|_{\Psi} := \inf \{ L>0 : \EE \Psi (Z/L)  < 1 \} .$$

We will assume that the noise is sub-Gaussian. Extension to sub-exponential noise
is straightforward but omitted to avoid technical digressions.

\begin{condition} \label{subGaussian.condition}
The error $\epsilon_1 $ is
independent of $ (X_1,Z_1)$ and satisfies for some constant $K_{\epsilon}\ge 1$
$$ \| \epsilon_1 \|_{\Psi} \le K_{\epsilon} .$$
\end{condition}

Recall that $P$ denotes the distribution of $(X,Z)$.
Let $p:= dP / d \nu $ be the density of $P$ with respect to a dominating product measure
$\nu := \nu_1 \times \nu_2$
with marginal densities $p_1$ and $p_2$. 
We define
$$r(x,z) := { p(x,z) \over p_1(x) p_2(z) } . $$

We let
$$ \gamma^2 := \int (r-1)^2 p_1 p_2 d \nu  $$
(assumed to exist). Note that $\gamma$ is the $\chi^2$-``distance" between the densities $p$ and 
$p_1 p_2$.

We impose the following incoherence condition.

\begin{condition} \label{incoherence.condition}
It holds that $\gamma < 1 $.
\end{condition}

Define
 $$f_{\rm P} := E ( f(X_1) \vert Z_1= \cdot ) , \ f_{\rm A} := f - f_{\rm P} . $$
 The subscript ``${\rm P}$" stands for ``projection", and
 ``${\rm A}$" stands for ``anti-projection". 
 Note that $f_{\rm P}$ is a function with the support of $Z_1$ as domain. 
 We assume this function to be smooth.
 
  \begin{condition} \label{Gamma.condition}
 For some constant $\Gamma$ it holds that
 $$J(f_{\rm P}) \le \Gamma \| f \| . $$
\end{condition}

To illustrate this condition, suppose 
that $Z_1$ is real-valued and 
$J(g) = \int | g^{(m)} (z) | dz $.
Suppose moreover that$$\sup_{x} \int | p^{(m)} (x \vert z) | dz \le \Gamma ,$$
where $p^{(m)} (x \vert z) := {d^m \over dz^m } ( p(x,z) / p_2(z)) $.
Then, interchanging differentiation and integration (and assuming this is allowed)
$$J(f_{\rm P} ) = \int \biggl  |  \int f(x) { p^{(m)} (x\vert z) } d \nu_1 (x) \biggr  |  dz  \le  \Gamma \int |f(x) | 
d \nu_1 (x) \le \| f \| .  $$

\section{Main result} \label{main.section}

We define 
\begin{equation}\label{tau.equation}
\tau_R(f,g) := \| f+g \| + \lambda I (f) + (\mu/R)^{2- q \over q} \mu J (g) . 
\end{equation}
We moreover let
\begin{equation} \label{tauI.equation}
\tau_I^2 (f) := \| f \|^2 + \lambda^2 I^2 (f) . 
\end{equation} 
 
 \begin{theorem} \label{main.theorem} 
 Assume Conditions  \ref{entropy.condition}, \ref{supnorm.condition}, 
\ref{subGaussian.condition}, \ref{incoherence.condition} and \ref{Gamma.condition}.
 Suppose that for some
 $0< \delta  < 1 $, $\max \{ A_I^2, A_J^2 \}/ n   \le n^{- \delta } $ and
 $(A_I^2 /  n)^{1 \over 1+ \alpha} \le (A_J^2 / n )^{1 \over 1+ \beta} n^{-\delta } $. 
 There exist  a universal constant $C$ and 
 constants $c$, $c_0$, $c_1$, $c_2$ depending on $\alpha$, $\beta$, $\gamma$, $\delta$, $B$, $\Gamma$, $q$
 and $K_{\epsilon}$ as well as on $I(f^0)$ and $J(g^0) $
such that for  $n \ge c_0$ and
 $$ \sqrt n \lambda^{1+ \alpha} = c_1 A_I , \ \sqrt n \mu^{1+ \beta} = c_1  A_J , $$
 $$ R  =  c_2  \mu   , \ R_I =  c_2  \lambda  , $$
 one has
$$\PP \biggl ( \tau_R ( \hat f - f^0 , \hat g - g^0 ) \le R , \ \tau_I ( \hat f - f^0 ) \le R_I \biggr ) \ge 1- C \exp[-n 
\lambda^2 / c ] . $$
\end{theorem} 

The proof is given in Section \ref{proofs.section}.

Theorem \ref{main.theorem}  does not provide the explicit dependence on
the constants. This dependence can in principle be deduced from Lemmas \ref{TR.lemma} and \ref{TRI.lemma} albeit that the expressions are somewhat complicated.
In an asymptotic formulation, considering $\alpha $, $\beta $, $\gamma$, $\delta$, $B$, $\Gamma$, $q$,
$K_{\epsilon}$ as well as $I(f^0)$ and  $J(g^0)$, 
as fixed, we get for $\lambda^2 \asymp A_I^{2 \over 1+ \alpha} n^{-{1 \over 1+\alpha}} $ and $\mu^2 \asymp 
 A_J^{2 \over 1+ \beta} n^{-{1 \over 1+\beta }}$, the rates
$$ \| \hat f - f^0 \|^2 = {\mathcal O}_{\PP} ( A_I^{2 \over 1+ \alpha}  n^{-{1 \over 1+ \alpha}} ) ,\  \| \hat g - g^0 \|^2 = {\mathcal O}_{\PP} ( A_J^{2 \over 1+ \beta}  n^{-{1 \over 1+ \beta  }} ) , $$
$$ I (\hat f) = {\mathcal O}_{\PP} ( 1) , \ J(\hat g) = {\mathcal O}_{\PP} (1) . $$

\begin{example} \label{Sobolev.example}
Suppose that $X_1$ and $Z_1$ take values in the interval $[0,1]$ and that
$I^2 (f) = \int | f^{(k)} (x) |^2 dx $ and $J^2 (g) = \int | g^{(m)} (z) |^2 dz $ with $m < k$. Then with $q=2$
the estimator is a spline and easy to calculate as the loss function as well
as the penalties are quadratic forms. The rates of convergence are
$\| \hat f - f^0 \|= {\mathcal O}_{\PP} (n^{-{k \over 2k+1}}) $ and
$\| \hat g - g^0 \|= {\mathcal O}_{\PP} (n^{- {m \over 2m+1 } })$.
See Section \ref{simulation.section} for some numerical results. 

\end{example} 

\begin{example} Suppose that $X_1$ takes its values in $[0,1]$ and $Z_1$ is real-valued.
Let $I^2 (f) :=  \int | f^{(k)} (x) |^2 dx $ with $k > 1$ and $J(g) := {\rm TV} (g)$ be the total variation of $g$.
Then with $q=1$ the estimator is again easy to calculate (the problem being formally equivalent to
a Lasso problem). The rates of convergence are
$\| \hat f - f^0 \|= {\mathcal O}_{\PP} (n^{-{k \over 2k+1}}) $ and 
$\| \hat g - g^0 \|= {\mathcal O}_{\PP} (n^{- {1 \over 3 } }\log^{1 \over 3}  n)$.
Indeed, Condition \ref{entropy.condition} for the class ${\cal G}$ now holds with
$\beta = 1/2$ and $A_J \asymp \sqrt {\log n }$.  This follows e.g. from
Lemma 2.2 in \cite{vandeGeer:00}. We note that once we have this fast rate for
$\| \hat f - f^0 \|$, the $(\log n )$-term in the rate for $\| \hat g - g^0 \|$ can be easily removed using
instead of the uniform entropy ${\cal H}_n$
the $\| \cdot \|_n$-entropy bound from  \cite{Birman:67} with $\| \cdot \|_n$-being the
empirical $L_2$-norm (i.e. for a real-valued function $m$ on
the support of  $(X_1, Z_1)$, $\| m \|_n := ( \sum_{i=1}^n m^2 (X_i, Z_i) / n )^{1/2}$). 
\end{example}

\section{Simulation results}\label{simulation.section}

\begin{figure}[H] 
\centering 
\begin{subfigure}[b]{0.45\textwidth}
 \centering
\includegraphics[width=\textwidth]
{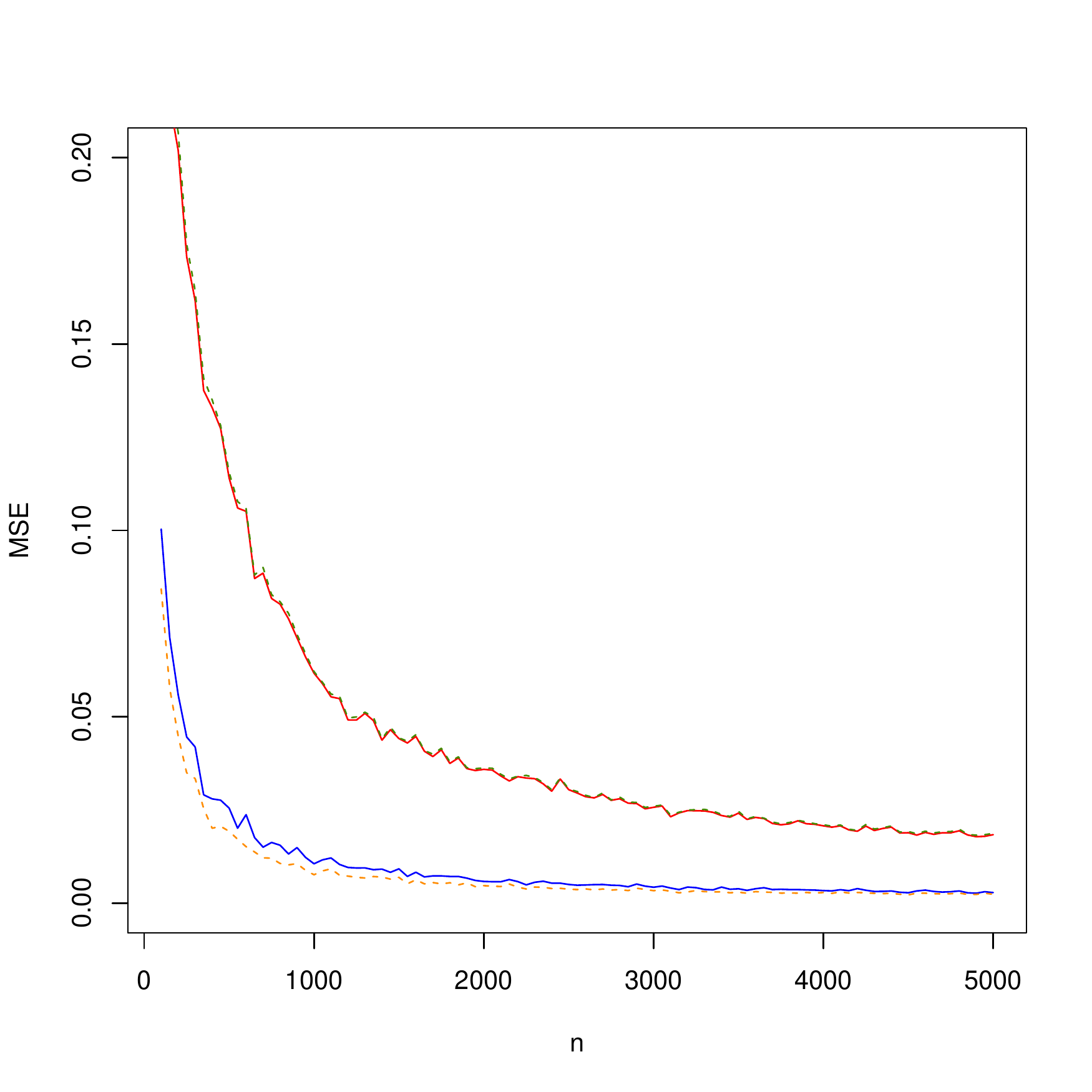}

\caption{$\rho\,=\,0.2$, SNR = 0.5}
\end{subfigure} 
~ 
\begin{subfigure}[b]{0.45\textwidth} \centering
\includegraphics[width=\textwidth]
{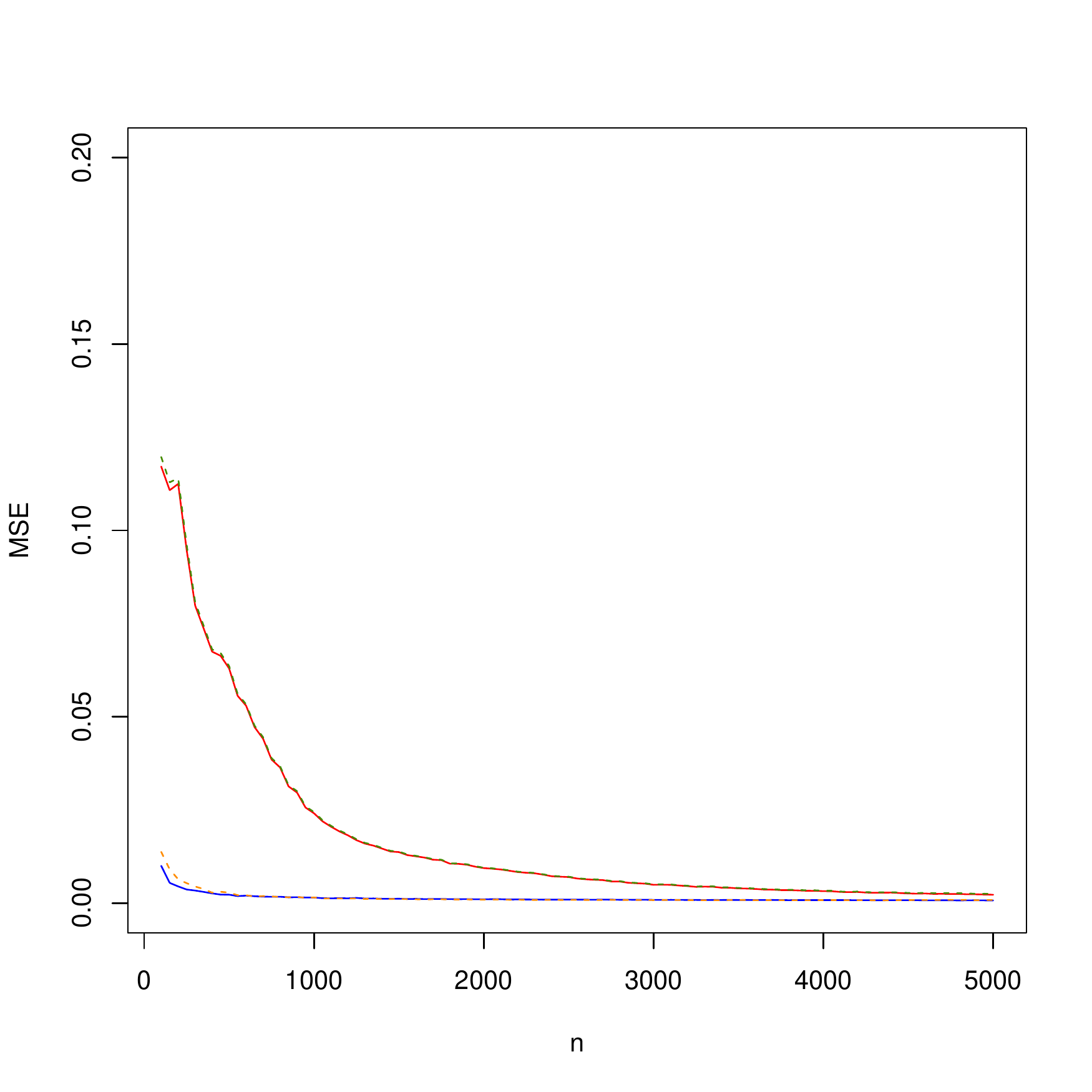}
\caption{$\rho\,=\,0.2$, SNR = 7} \label{subfig:bestMISE} 
\end{subfigure} 

\begin{subfigure}[b]{0.45\textwidth}
\centering 
\includegraphics[width=\textwidth]
{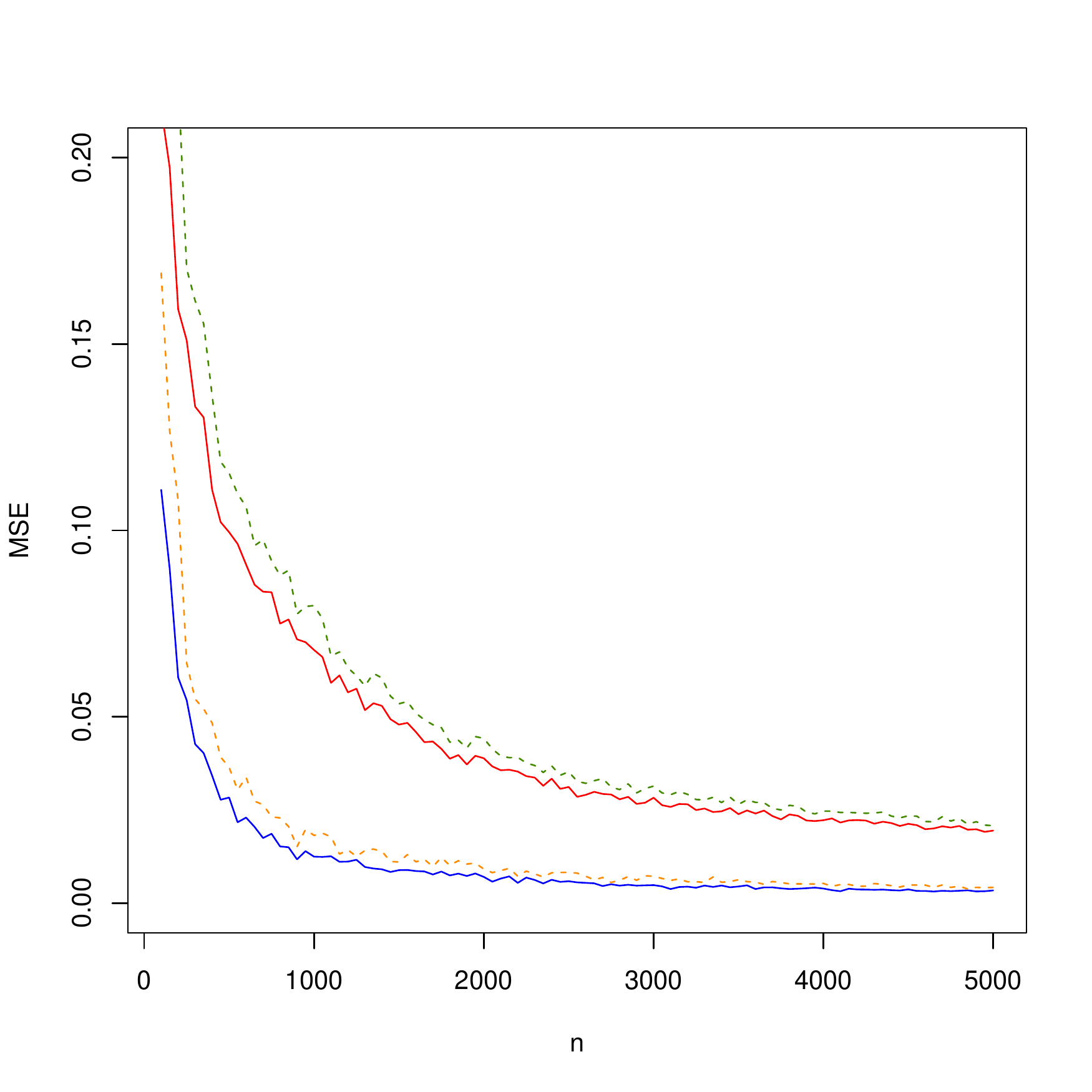}
\caption{$\rho\,=\,0.8$, SNR = 0.5}
\end{subfigure} 
~ 
\begin{subfigure}[b]{0.45\textwidth} 
\centering
\includegraphics[width=\textwidth]
{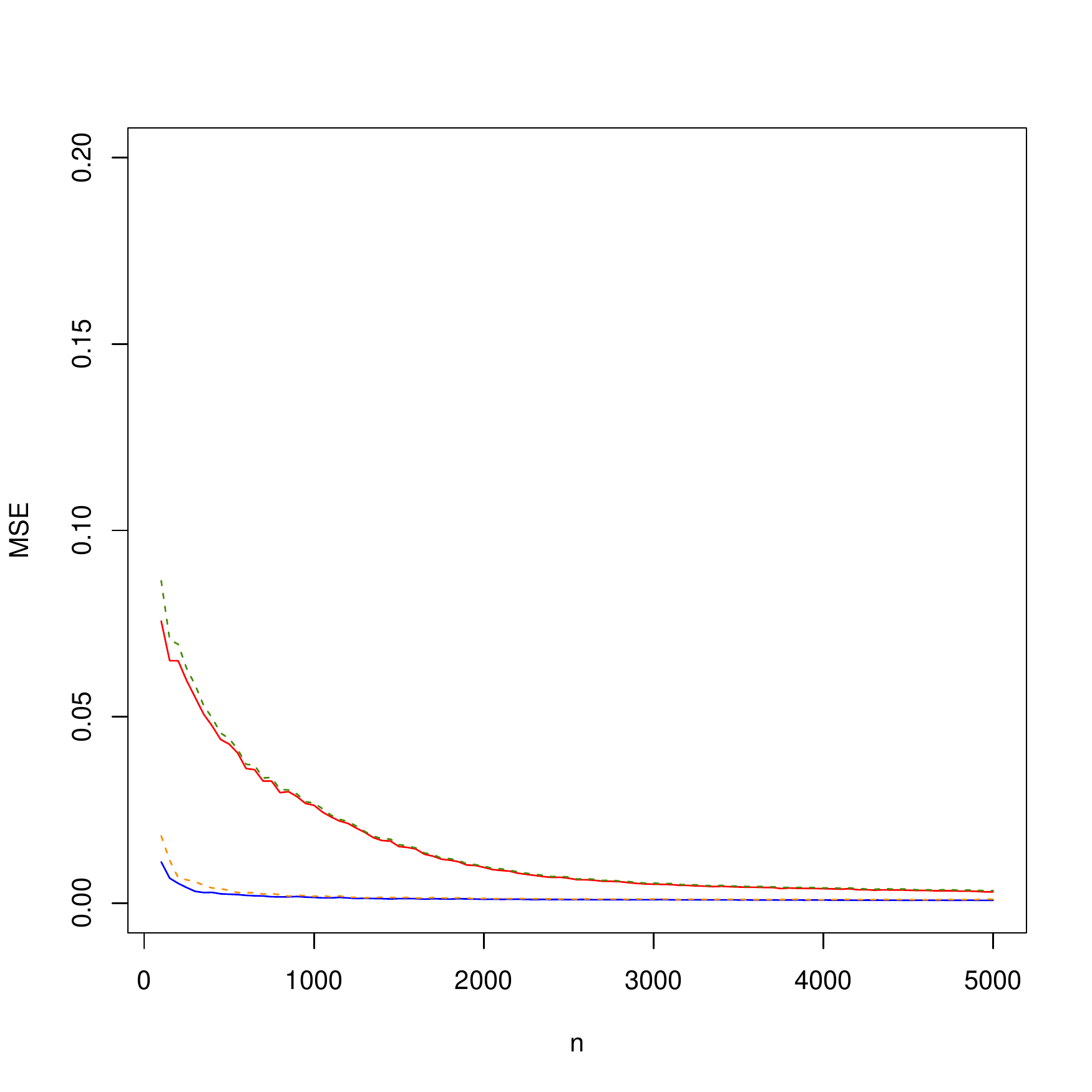}
\caption{$\rho\,=\,0.8$, SNR=7}
 \end{subfigure}
\caption{Estimated MSE for each of the computed estimators: $\hat{f}^s$ (blue line), $\hat{g}^s$ (red line), $\hat{f}$ (orange dotted line), and $\hat{g}$ (green dotted line) for the four simulation scenarios.}\label{fig:MSE} 
\end{figure}

In this simulation study, we show that the results of Theorem \ref{main.theorem} also (approximately) hold empirically. We consider Example \ref{Sobolev.example}. We estimate each of the ``true'' functions $f^0$ and $g^0$ in the cases where neither functions are known and the cases where one of them is known. We will see that, for each function, the rate of convergence of the estimator when neither of the ``true" functions is known is of the same order than that when one of the components is known. For this, we will show the plots of the MSE of the four estimators in four different scenarios (see Figure \ref{fig:MSE}). However, we will only show the plots of the estimators when correlation$(X,Z)\,=\,0.8$, SNR $=\,7$ since analogous results hold for the other scenarios.

Let $X$ and $U$ be independent uniformly
distributed random variables with values in $(0,1)$. Define $Z\,=\,a\,X+ (1-a)\,U$ with $a$ an
appropriate constant such that the correlation between X and Z is equal to
$\rho$ (which we will define later). 

We use B-splines of order 6 (piecewise polynomials
of degree 5) to represent each of the functions $f$ and $g$
(see \cite{depractical}). We write

\begin{equation*} 
f(x) = \sum_{i=1}^{K} \gamma_{f,i} b_{f,i} (x),\quad\quad g(z) =
\sum_{j=1}^{K} \gamma_{g,j} b_{g,j} (z),
 \end{equation*} 
~\\
 where $b_{f,i},b_{g,j},\,\,i,j=1,...,K$ are the basis functions of the
 B-spline parametrization, $\gamma_{f}=(\gamma_{f,1},...,\gamma_{f,K}),
 \gamma_{g}=(\gamma_{g,1},...,\gamma_{g,K})$ are the parameters vectors of
 $f$ and $g$, respectively, and $K+6$ is the number of knots, which we choose
 to be $\lceil 3 \sqrt n /5 \rceil +6$ where $n$ represents the number of observations.
 Denote by
$(x_1,...,x_n)$ and $(z_1,...,z_n)$ realizations of the dependent random
variables $X$ and $Z$ and let
$x_{(r)}$ be  the $r$-th order statistic of the sample from $X$ ($r=1 , \ldots , n$). 
For estimating
 the function $f$ (and analogously for the function $g$), we place the first and last $6$ knots (corresponding to
 the order of the B-spline) in $x_{(1)}$ and
 $x_{(n)}$, respectively, and position the remaining knots uniformly in $\{x_{(2)},...,x_{(n-1)}\}$. 
%
%
We define the penalizations as
\begin{equation*} 
I^2(f) := \int_0^1 | f'''(x) |^2 \mathrm{d}x \; + \; \int_0^1 | f(x) |^2
\mathrm{d}x, \end{equation*}
\begin{equation*} 
J^2(g) := \int_0^1 |g''(z) |^2 \mathrm{d}z \; + \; \int_0^1 |g(z)
|^2 \mathrm{d}z
\end{equation*} 

and the $(i,j)-th$ components of the matrices $\Omega_f, \Omega_g\in\mathbb{R}^{K \times K}$ as

\begin{equation*}
 (\Omega_f)_{i,j} := \int_0^1 b_{f,i}'''(x)b_{f,j}'''(x) \mathrm{d}x +
\int_0^1 b_{f,i}(x)b_{f,j}(x) \mathrm{d}x 
\end{equation*}
 and
\begin{equation*} 
(\Omega_g)_{i,j} := \int_0^1 b_{g,i}''(z)b_{g,j}''(z) \mathrm{d}z +
\int_0^1 b_{g,i}(z)b_{g,j}(z) \mathrm{d}z
\end{equation*} 
~\\
Then, we
can write $I^2(f) = \gamma_f^T \Omega_f \gamma_f$ and $J^2(g) = \gamma_g^T \Omega_g
\gamma_g$. Moreover, using Cholesky, we
can find matrices $H_f, H_g \in \mathbb{R}^{K \times K}$ such that $\Omega_f = H_f^T H_f$
and $\Omega_g = H_g^T H_g$. 
~\\
 ~\\
 \underline{\textbf{The case where both $f^0$ and $g^0$ are unknown:}}
~\\
Consider the two-components model:

\begin{equation}\label{eq:Original} 
Y_i = f^0 (X_i) + g^0 (Z_i) + \epsilon_i, \quad\quad i=1,...,n,
\end{equation} 
~\\
where $\epsilon_i,\;i=1,...,n$ are i.i.d. centered Gaussian random variables with variance $\sigma^2$. 
The estimator is
$$(\hat f , \hat g) := \arg \min_{f,g} \biggl \{ 
\| Y - f - g \|_n^2 + \lambda^2 I^2 (f) + \mu^2 J^2 (g) \biggr \} . $$

%
We took $\lambda = 14\,n^{-3/7}$ and $\mu = 0.3\,n^{-2/5}$.  The
constants of both tuning parameters are chosen by minimizing the mean square 
error\footnote{Estimated using 100 simulations.} of the estimators for the case $n=5000$. Candidates for the constants were taken from the grid
$(\{1,2,3,...,20\}\times\{0.1,0.2,0.3,...,1\})$, where the first set
corresponds to the constant of $\lambda$ and the second to the constant of
$\mu$.

%
%
\underline{\textbf{The case where $f^0$ or $g^0$ is known:}}
~\\ 
If $g^0$ is known we re-write equation \eqref{eq:Original} as
~\\
\begin{equation*}
 Y^f_i = f^0 (X_i) + \epsilon_{i}, 
\end{equation*}
with $ Y^f_i = Y_i - g^0(Z_i)$, $i=1 , \ldots , n $. 

%
We then use the estimator
$$ \hat{f}^{s} := \arg \min_f \biggl \{  \| Y^f - f \|_n^2 + \lambda^2 I^2 (f)\biggr \}  .$$
The tuning parameter is taken to be $\lambda = 14\,n^{-3/7}$. 

Similarly, if $f^0$ is known we let $Y^g:= Y-f^0$ and
$$\hat{g}^{s} := \arg \min_g \biggl \{  \| Y^g - g \|_n^2 + \mu^2 J^2 (g)  \biggr \} $$
with $\mu =0.3\,n^{-2/5}$.

\textbf{Simulations}:  
~\\ 
Define the Signal-to-Noise ratio as ${\rm SNR} :={\rm var}(f^0(X) + g^0(Z))/{\sigma^2}$. For our simulations, we consider the following scenarios:

\begin{itemize} 
\item $f^0 (x) = -10 \sin (1.9x + 0.2\pi) + \mathbb{E}[10 \sin (1.9x +
0.2\pi)]$. 
\item $g^0(z) = 3 e^{ -500(z - 0.1)^2 } - \mathbb{E}[3e^{ -500(z - 0.1)^2 }]$.
\item SNR $\in\{ 0.5,7 \}$. 
\item $\rho\,\in\,\{0.2,0.8\}$\footnote{The value $\rho=0.2$ corresponds to $a=0.169$
  and the value $\rho=0.8$ to $a=0.571$.}. 
\item $n\in \{100,150,200,...,5000\}$. 
\end{itemize} 

The error variance $\sigma^2$ was chosen in each scenario to match the above given
Signal-to -Noise ratios. 
For each $n$ the average of 100 simulations is used to estimate the mean
square error.
In Figure \ref{fig:MSE}, we see that the rate of convergence of $\hat{f}$
and of $\hat{f}^s$ are of similar order and that the same applies to
$\hat{g}^s$ and $\hat{g}$. In other words, for each function $f^0$ and
$g^0$, the rate of convergence of the estimators when both functions are
unknown (approximately) corresponds to the case when one of them is
known. These results agree with Theorem \ref{main.theorem} and hold in
the four simulation scenarios. Moreover, we see that the convergence of
$\hat{f}^s$ and $\hat{f}$ to $f^0$ is faster than that of $\hat{g}^s$ and
$\hat{g}$ to $g^0$, which is also established in Theorem \ref{main.theorem}. 

\begin{figure}[H] 
\centering
\begin{subfigure}[b]{0.45\textwidth} 
\centering 
\includegraphics[width=\textwidth]
{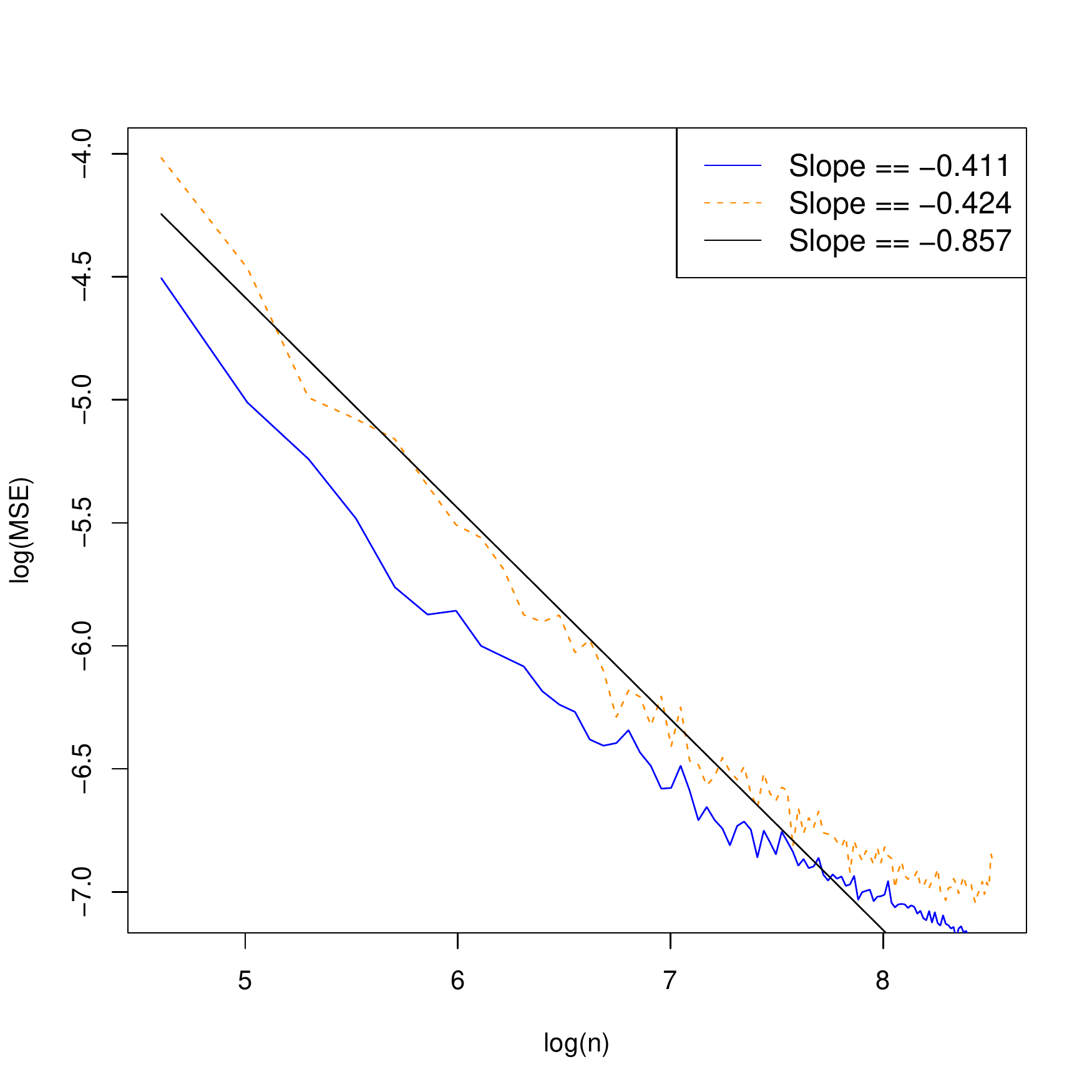}
 \caption{Log transformed data for $\hat{f}^s$ (blue line) and $\hat{f}$ (orange dotted line)}
\end{subfigure} 
~ 
\begin{subfigure}[b]{0.45\textwidth} 
\centering 
\includegraphics[width=\textwidth]
{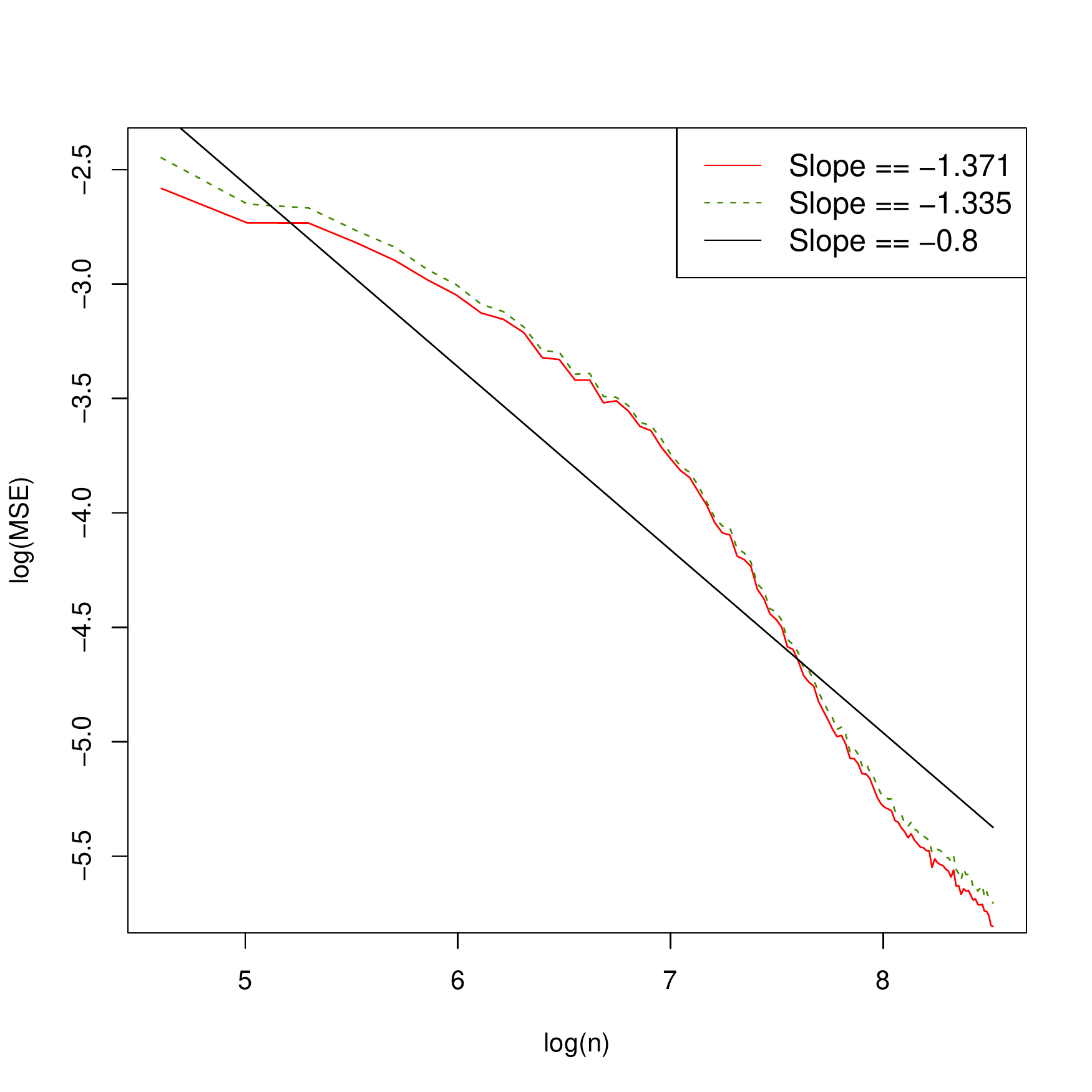}
 \caption{Log transformed data for $\hat{g}^s$ (red line) and $\hat{g}$  (green dotted line)}
\end{subfigure} 
\caption{Log-transformed data for the case $\rho\,=\,0.8$ and SNR = 7. A black line using the theoretical slope and an arbitrary intercept was drawn for graphical comparison.}\label{fig:Corr0.8SNR7rd}
\end{figure}

The log-transformed data from Figure \ref{fig:MSE} for the scenario $\rho\,=\,0.8$ and SNR = 7 is plotted in Figure
\ref{fig:Corr0.8SNR7rd}. Here, we fit a linear regression on each curve considering only 
those observations corresponding to $n\,\geq\,1000$ and print the slope of
these and the theoretical slope\footnote{Recall that by Theorem \ref{main.theorem}
  we have $\log ||\hat{f} - f^0 ||^2_2 =\log(c_1) -(6/7)\log(n)$ and $\log
  ||\hat{g} - g^0 ||^2_2 =\log(c_2) -(4/5)\log(n)$, where $c_1$ and $c_2$
  are constants depending on those of the tuning parameters.} in
the legend of the plot. With SNR=7 it is not clear whether the slopes of the regression line of the estimators  
 agree with their theoretical counterpart. For lower SNR however the agreement is remarkably
good (not shown here). 

\begin{figure}[H] 
\centering 
\begin{subfigure}[b]{0.45\textwidth} 
\centering 
\includegraphics[width=\textwidth]
{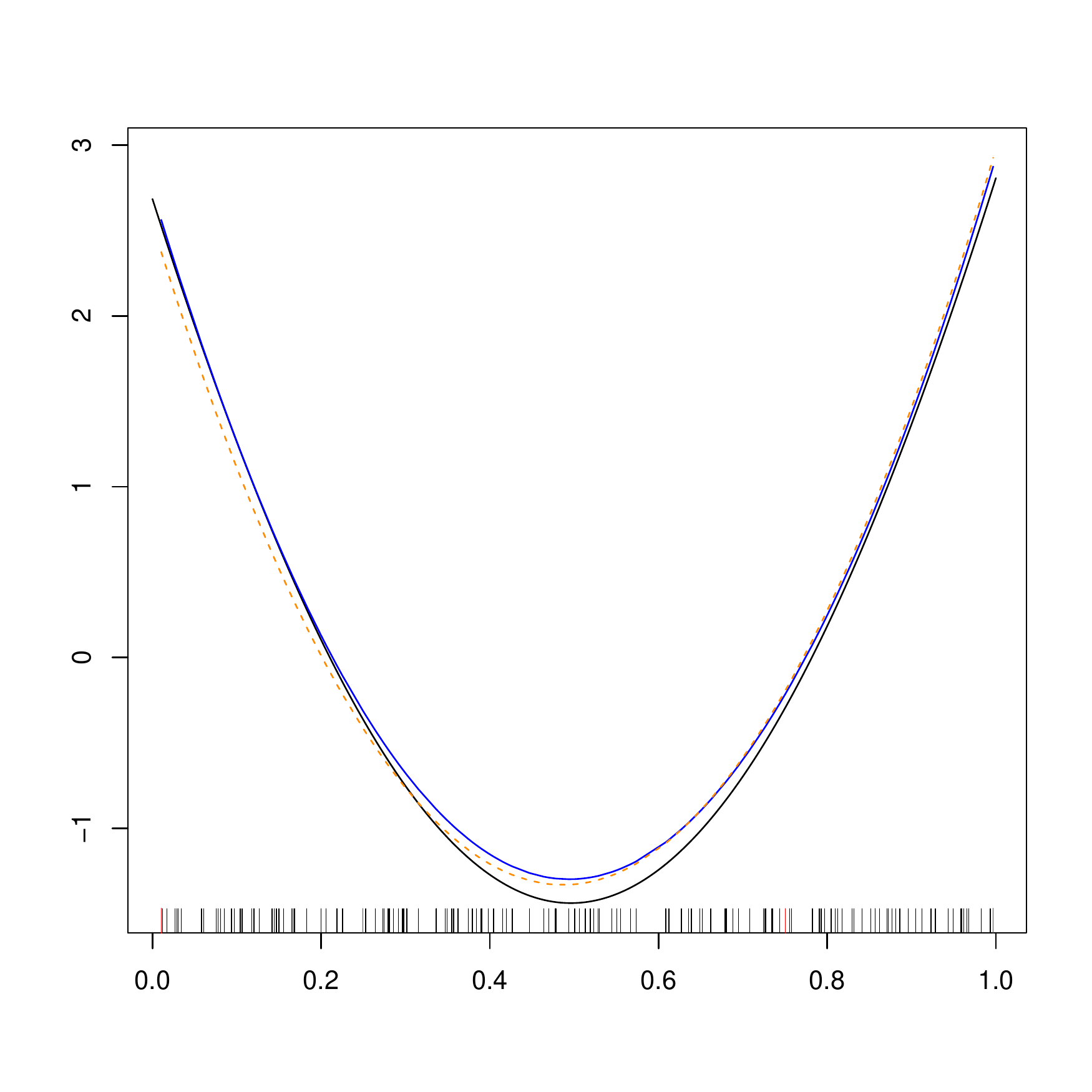}
\caption{$n$ = 150}
\end{subfigure}
 ~ 
\begin{subfigure}[b]{0.45\textwidth} 
\centering 
\includegraphics[width=\textwidth]
{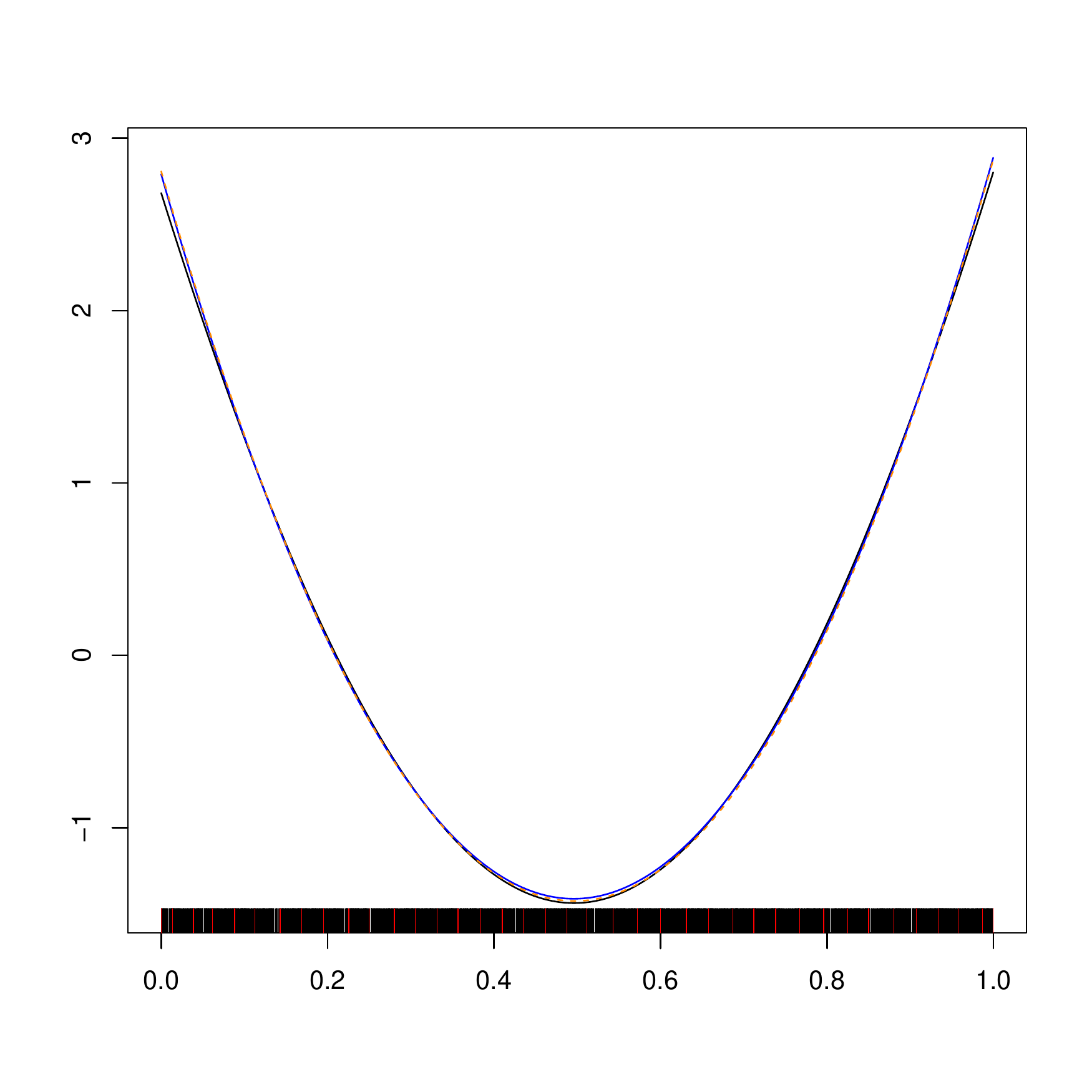}
 \caption{$n$ = 5000}
\end{subfigure} 
\\ 
\begin{subfigure}[b]{0.45\textwidth} 
\centering 
\includegraphics[width=\textwidth]
{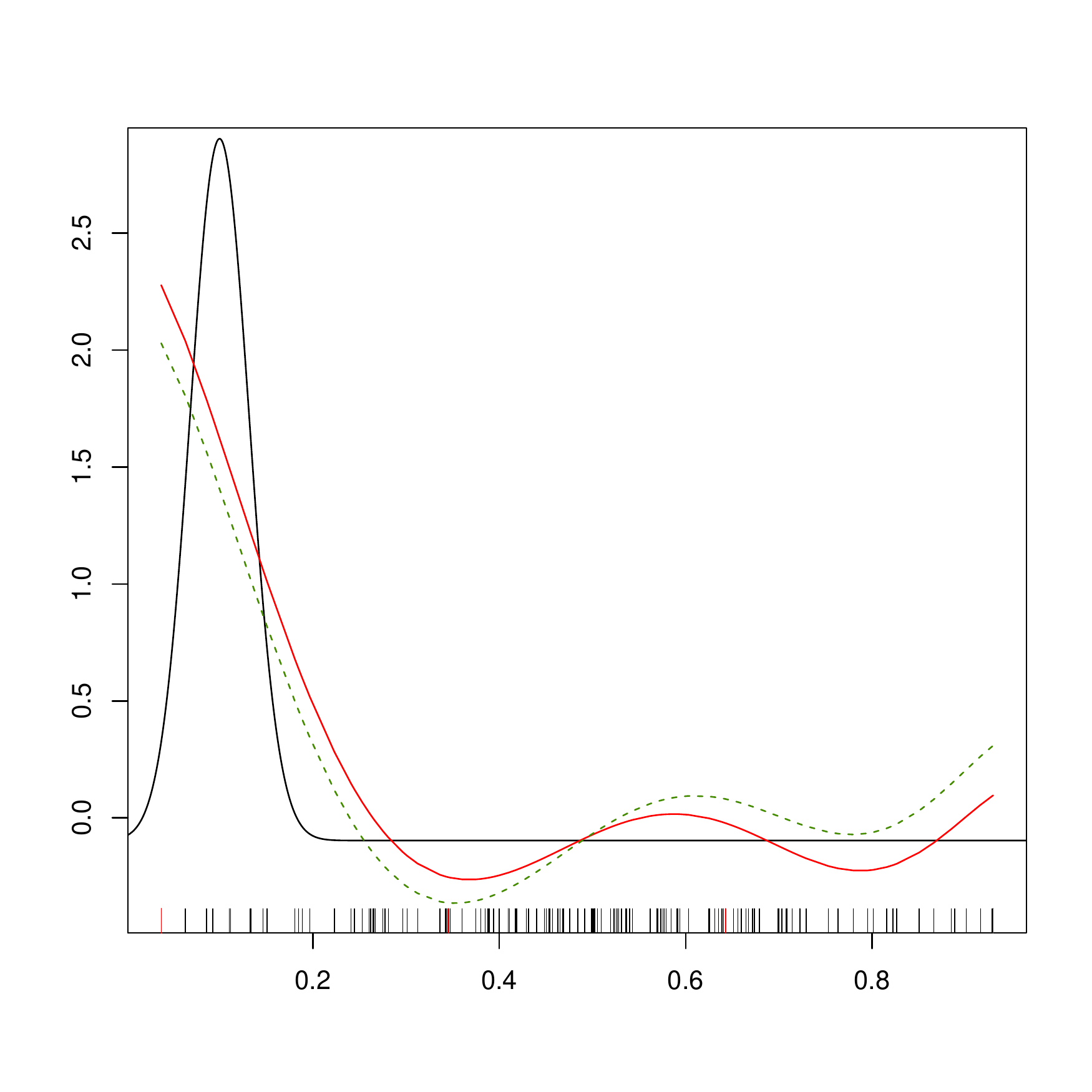}
 \caption{$n$ = 150}
\end{subfigure} 
~ 
\begin{subfigure}[b]{0.45\textwidth} 
\centering 
\includegraphics[width=\textwidth]
{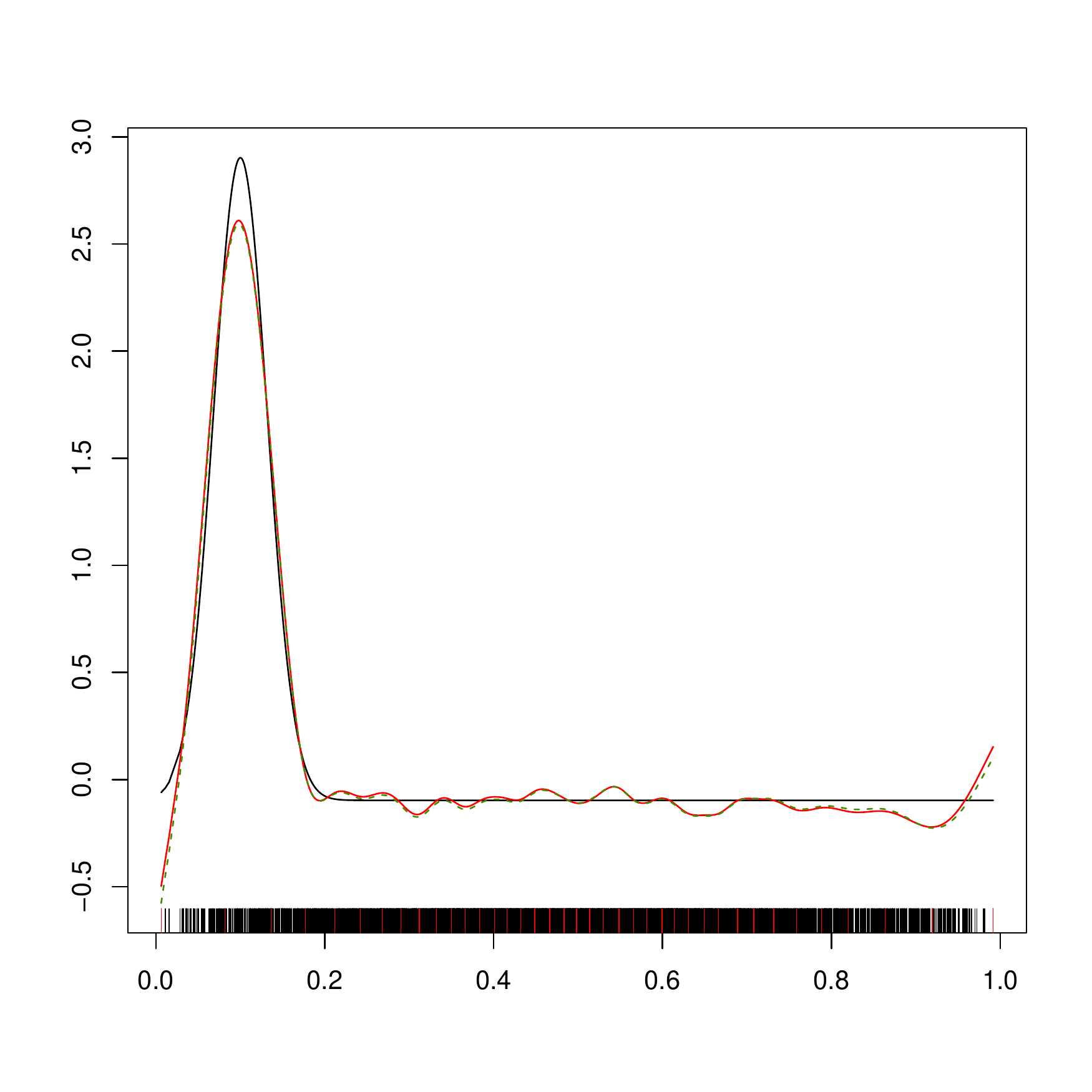}
 \caption{$n$ = 5000}
\end{subfigure} 

\caption{Plots of the true functions $f^0$ and $g^0$ (black lines) with their corresponding estimators  $\hat{f}^s$ (blue lines), $\hat{f}$ (orange dotted lines), $\hat{g}^s$ (red lines), $\hat{g}$ (green dotted lines), for $\rho\,=\,0.8$ and SNR = 7. The data are represented
with black small vertical lines and knots positions with red small vertical
lines. For each $n \in \{150,5000\}$, we use a single simulation (not 100 simulations).}
\label{fig:Corr0.8SNR7_fg} 
\end{figure}

\begin{figure}[H]
\centering
\begin{subfigure}[b]{0.45\textwidth} 
\centering 
\includegraphics[width=\textwidth]
{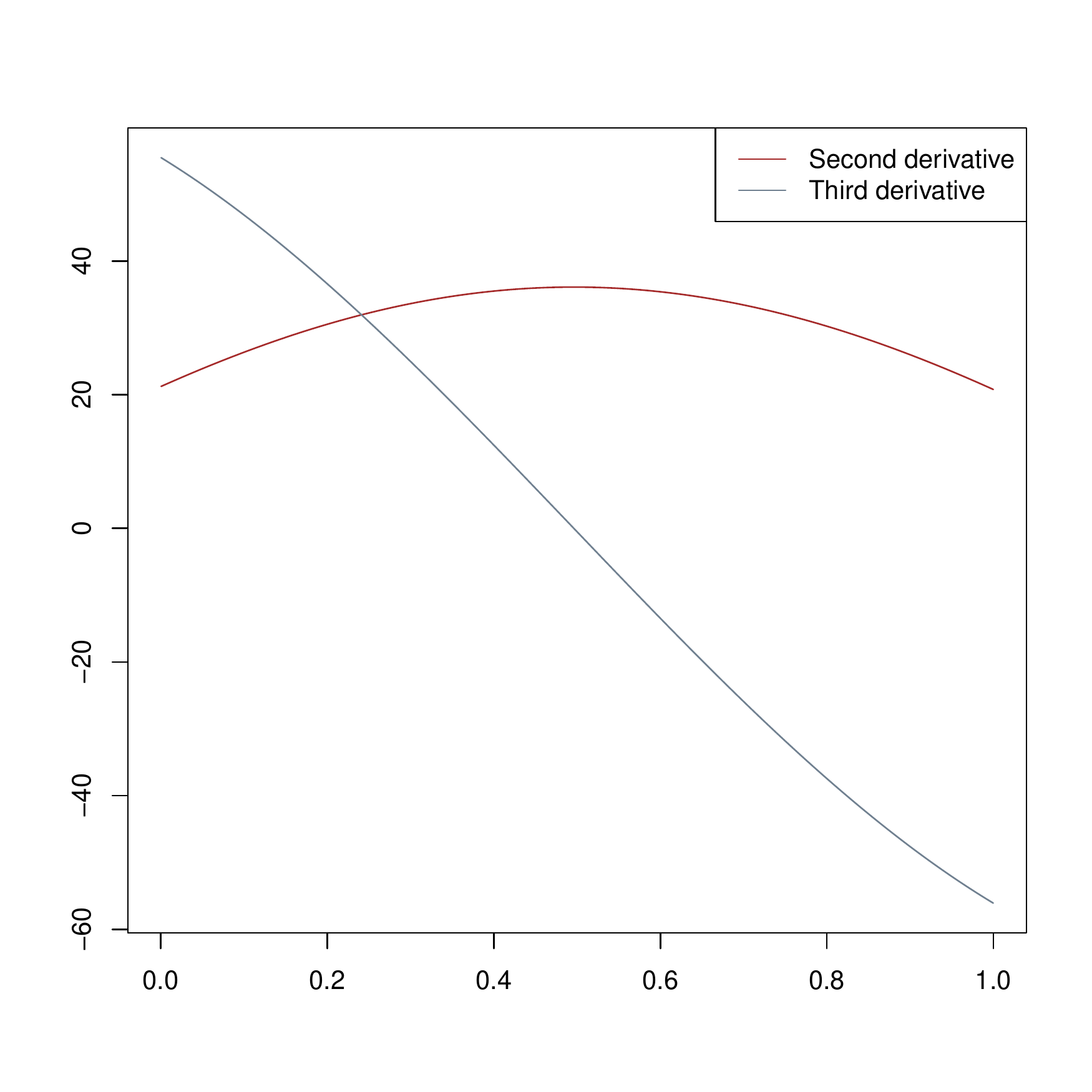}
 \caption{Derivatives of $f^0$.}
 \end{subfigure} 

\begin{subfigure}[b]{0.45\textwidth} 
\centering 
\includegraphics[width=\textwidth]
{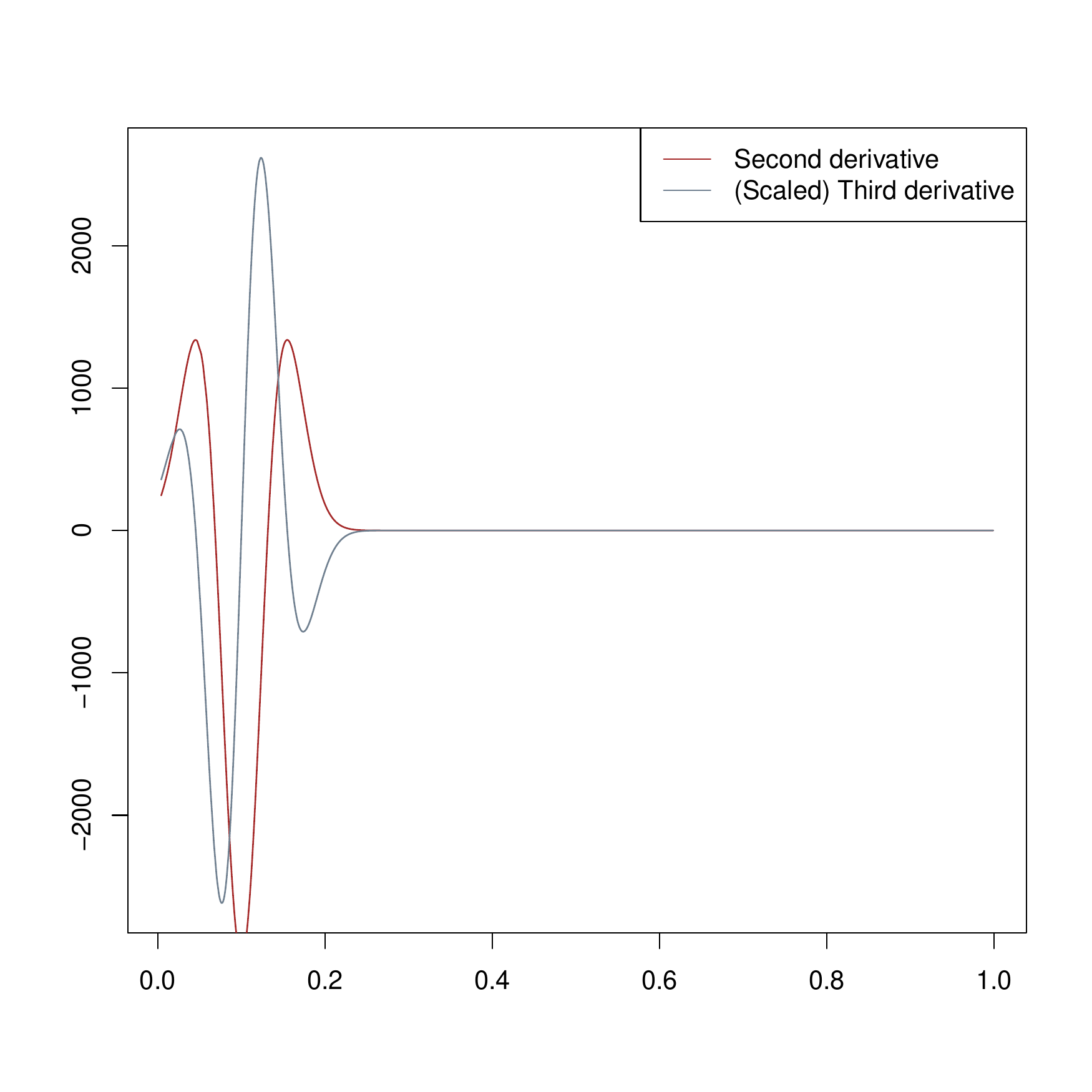}
 \caption{Derivatives of $g^0$.}
 \end{subfigure} 
\caption{ Second and third derivatives of the true functions $f^0$ and $g^0$ in $[0,1]$. The third derivative of $g^0$ was multiplied by 0.02 for an easier comparison.} \label{fig:derg0} 
\end{figure}

The plots of both $f^0$ and $g^0$ and their corresponding estimators for
the scenario $\rho\,=\,0.8$ and SNR = 7 are displayed in Figure
\ref{fig:Corr0.8SNR7_fg}. We can see that, as the number of observations
increases, the functions $\hat{f}$ and $\hat{g}$  converge to $\hat{f}^s$
and $\hat{g}^s$, respectively. This happens while all of them improve their
estimation of the true functions $f^0$ and $g^0$ appropriately. We note
that $\hat{f}$ and $\hat{f}^s$ are almost identical to $f^0$ when the
number of observations is large. However, $\hat{g}$ and $\hat{g}^s$ can
only resemble but not describe perfectly $g^0$. This is probably due to the highly
variable second and third derivatives of $g^0$ in comparison with those of $f^0$, as can be seen in Figure \ref{fig:derg0}.

 ~\\ 

%
%
%

%
\section{Proofs} \label{proofs.section}

We use the notation
$P_n := \sum_{i=1}^n \delta_{(X_i, Z_i)} / n $
for the empirical measure based on $\{ (X_i , Z_i ) \}_{i=1}^n $.

The proof is organized as follows. We first present some preliminary results needed for
the proof of the faster rate for $\hat f$. Then we look in Subsection \ref{global.section} at the global rate for
both components. We use here the convexity of the least squares loss function
and the penalties to localize the problem to the set ${\cal M} (R) :=
\{ (f,g) :\ \tau_R (f,g ) \le R \}$, and then show that indeed $(\hat f - f^0 , \hat g - g^0) \in {\cal M} (R)$
provided that the random part of the problem is under control. In Subsection  \ref{TR.section}
we show the random part is indeed under control with large probability. For this result,
we need recent findings from empirical process theory, in particular the convergence
of empirical norms and inner products. Here, we apply some results from
\cite{vdG2013causal}. The application is somewhat elaborate: for an
additive model with $p$ components there are ${ p+1 \choose 2 } -1$ terms to consider.
If there is only one component, say $f$, one needs to consider the behaviour
of $\| f\|_n^2 - \| f \|^2 $ and $\epsilon^T f / n $ uniformly over some collection of functions
$f$. If there are two components $f$ and $g$ the number of terms to consider is five:
namely uniform convergence
$\| f\|_n^2$, $\| g\|_n^2$, $P_nfg$, $\epsilon^T f/n $ and $\epsilon^T g /n$ to their theoretical
counterparts. This is done in Subsection
\ref{empirical.section}. Subsection \ref{global.section} takes such uniform convergence for granted.
The same is true in Subsection \ref{tighter.section} where we show the faster rate
for the estimator $\hat f$ of the smoother component: the results are on a random event which
is shown to have large probability in Subsection \ref{TRII.section}
using results from
empirical process theory given in Subsection \ref{empirical.section}. We finally collect all
pieces in Subsection \ref{finish.section} to finish the proof of the main result in Theorem \ref{main.theorem}.

\subsection{Preliminaries}\label{preliminaries.section}

%
%

\begin{lemma} \label{incoherence.lemma}
Assume Condition \ref{incoherence.condition} and 
suppose $\int f p_1 d \nu_1   = 0$.
Then
$$ \| f + g \|^2 \ge  (1-\gamma) ( \| f \|  + \| g \|)^2 . $$
\end{lemma}

{\bf Proof.}
We have
$$ \| f + g \|^2 = \| f \|^2 + \| g \|^2 + 2 \int (f g) p d \nu  . $$
Moreover, since $\int f p_1  = 0$, 
$$| \int fg p d \nu  |= | \int fg  (r-1) p_1 p_2 d \nu  | \le \gamma ( \int f^2 g^2 p_1 p_2 d \nu  )^{1/2} = \gamma \| f \| \| g \| . $$ 
Hence,
$$ \| f + g \|^2 \ge \| f \|^2 + \| g \|^2- 2\gamma \| f \| \|g \| $$
$$ = (1- \gamma) (  \| f \|^2 + \| g \|^2) + \gamma ( \| f \| - \| g \| )^2  \ge (1- \gamma) ( \| f \|  + \| g \|)^2 . $$
\hfill $\sqcup \mkern -12mu \sqcap$

 \begin{lemma} \label{anti.lemma} Assume Condition \ref{incoherence.condition} and 
suppose $\int f p_1 d \nu_1 =0$. We have that
 $$ \| f_{\rm P} \|^2 \le \gamma \| f \|^2 $$ and
 $$\| f_{\rm A} \|^2 \ge (1- \gamma^2 ) \| f \|^2 . $$
 
 \end{lemma} 
 
 {\bf Proof.} We have
 $$ f_{\rm P}  = \int f (r-1) p_1 d \nu_1 .$$
 Hence
 $$\| f_{\rm P} \| \le  \| r-1 \| \| f \|  = \gamma \| f \| , $$
 and
 $$ \| f_{\rm A} \|^2 = \| f \|^2 - \| f_{\rm P}  \|^2 \ge (1- \gamma^2 ) \| f \|^2 . $$
 
 \hfill $\sqcup \mkern -12mu \sqcap$.
 
 \begin{lemma}\label{entropy.lemma} Assume Conditions \ref{entropy.condition} and
\ref{supnorm.condition}. Then
$${\cal J}_{\infty} (z , \{f_{\rm A} : \ f \in {\cal F} (R,M) \} ) \le 2 {\cal J}_{\infty} (z, {\cal F} (R,M) ), \ z > 0  $$
and for $R \le M / B$
$$\sup_{f \in {\cal F}(R,M) } \| f_{\rm A} \|_{\infty} \le 2  M . $$

\end{lemma}

{\bf Proof.} 
Let $u >0$ and $f  , \ \tilde f \in {\cal F} (R,M)$ be arbitrary, satisfying
$\| f - \tilde f \|_{\infty}  \le u$. 
Then clearly also
$$ \|  f_{\rm P} - \tilde f_{\rm P}  \|_{\infty} = \| E ( f (X_1) -\tilde f (X_1) \vert Z= \cdot )\|_{\infty} \le u. $$
So then
$$\| f_{\rm A} - \tilde f_{\rm A} \|_{\infty} \le \| f_{\rm P} - \tilde f_{\rm P}\|_{\infty} +  \| f- \tilde f  \|_{\infty}\le 2u . $$
Similarly, for $f \in {\cal F}(R,M)$, we have
$$\| f_{\rm A }\|_{\infty} \le \| f_{\rm P} \|_{\infty} +  \| f \|_{\infty} \le 2  M . $$
\hfill $\sqcup \mkern -12mu \sqcap$

 \subsection{A global bound} \label{global.section}
 
We define 
\begin{equation}\label{MR.equation}
{\cal M} (R):= \{ (f,g) : \ \tau_R (f,g) \le R \} 
\end{equation} and for a sufficiently small
value $\delta_0 >0 $, to be chosen later the sets
$${\cal T}_1 (R) := \biggl \{ \sup_{(f,g) \in {\cal M} (R) } \biggl | \| f+g \|_n^2 - \| f+g \|^2 \biggr | \le \delta_0^2 R^2 \biggr \} , $$
$${\cal T}_2 (R) := \biggl \{ \sup_{(f,g) \in {\cal M} (R) } | \epsilon^T (f+g)|/n  \le \delta_0^2 R^2  \biggr \}  ,$$
and
\begin{equation}\label{TR.equation}
{\cal T} (R) := {\cal T}_1 (R) \cap {\cal T}_2 (R) . 
\end{equation}

\begin{lemma} \label{global.lemma}Take $\delta_0 \le {1 \over 20}  $ and 
suppose that
\begin{equation} \label{IJbound.equation}
\lambda^2 I^2 (f^0) + \mu^2 J^q (g^0) \le  \delta_0^2 R^2 .
\end{equation}
Then on ${\cal T}(R)$, we have $\| \hat m - m^0 \|^2 + \lambda^2 I^2 (\hat f) +
\mu^2 J^q ( \hat g) \le 4 \delta_0^2 R^2 $ and $\tau_R ( \hat f - f^0 , \hat g - g^0)  \le R $.
\end{lemma}
{\bf Proof.}
Define
$$\tilde f := t \hat f + (1-t)f^0 , \ \tilde g := t \hat g + (1- t) g^0  $$
with
$$t:= { R \over R + \tau_R ( \hat f - f^0 , \hat g - g^0 ) } . $$
Then
$\tau_R( \tilde f  - f^0, \tilde g - g^0 ) \le R$. Let $\tilde m := \tilde f + \tilde g$ and
$m^0 := f^0 + g^0 $. By convexity 
$$\| \tilde m - m^0 \|_n^2 + \lambda^2 I^q ( \tilde f) + \mu^2 J^q (\tilde g) \le
2 \epsilon^T (\tilde m- m^0 ) + \lambda^2 I^q ( f^0) + \mu^2 J^q (g^0) . $$
On ${\cal T}(R)$ we find
$$\| \tilde m - m^0 \|^2 + \lambda^2 I^q (\tilde f) + \mu^2 J^q (\tilde g) \le
4 \delta_0^2 R^2  . $$
Hence
$$I ( \tilde f  ) \le ( 2 \delta_0 ) R / \lambda , $$
and
$$J( \tilde g) \le (2 \delta_0 )^{2/q} ( R/ \mu)^{2/q} \le 2 \delta_0 (R/\mu)^{2/q}$$
where in the last step we used $2 \delta_0 < 1$ and $2/q \ge 1 $. 
Since by (\ref{IJbound.equation}) it holds that $I(f^0)  \le 2 \delta_0 R/\lambda$ we get
$$I( \tilde f - f^0) \le 4 \delta_0 R/\lambda . $$
Also, by (\ref{IJbound.equation}) we have $J(g^0) \le (2\delta_0)^{2/q} ( R/\mu)^{2/q} \le 2 \delta_0 (R/\mu)^{2/q}$
so that
$$J( \tilde g - g^0 ) \le4 \delta_0 ( R/\mu)^{2/q}.$$
We find
$$\lambda  I ( \tilde f  - f^0 ) \le  4 \delta_0 R  $$
as well as
$$ (\mu/ R)^{2-q \over q }   \mu J( \tilde g - g^0 ) \le 4 \delta_0 R  . $$
But then
$$ \tau_R ( \tilde f  - f^0 , \tilde g - g^0)  = \| \tilde m - m^0 \| + \lambda  I (\tilde f-f^0) + (\mu/R)^{2-q \over q}  \mu J (\tilde g-g^0 ) $$ $$\le
10 \delta_0 R \le R/2 $$
where we used $\delta_0 \le {1 \over 20} $. 
This implies $\tau_R (\hat f - f^0 , \hat g - g^0 ) \le R $. Repeating
the argument completes the proof. 
\hfill $\sqcup \mkern -12mu \sqcap$

\subsection{A tighter bound for the smoother part} \label{tighter.section}

 Let ${\cal F}(R_I) := \biggl \{ f : \  \tau_I (f)  \le R_I  \biggr \}. $

 For $\delta_1 $ sufficiently small we define 
 $${\cal T}_{1,I} (R_I):= \biggl \{ \sup_{f \in {\cal F} (R_I)}  \biggl | \| f_A \|_n^2 - \| f_A \|^2 \biggr | \le \delta_1^2 R_I^2  \biggr \} ,$$
 $${\cal T}_{I,2} (R_I) := \biggl \{ \sup_{f \in {\cal F}(R_I)}  | \epsilon^T f_A    | /n \le \delta_1^2 R_I^2  \biggr \} , $$
 $${\cal T}_{I,3} (R_I, R):= \biggl \{ \sup_{(f,g)\in {\cal M} (R): \ f  \in {\cal F} (R_I)} |P_n f_A (g + f_P)    | \le \delta_1^2 R_I^2  \biggr \} $$
 and we let
 \begin{equation} \label{TRI.equation}
 {\cal T}_I (R_I, R) := {\cal T}_{I,1}(R_I)  \cap {\cal T}_{I,2} (R_I) \cap {\cal T}_{I,3} (R_I, R) .  
 \end{equation}
 
 \begin{lemma} \label{tighter.lemma} Assume Condition \ref{incoherence.condition} and
\ref{Gamma.condition}. Suppose the condition  (\ref{IJbound.equation})
$$\lambda^2 I^2 (f^0) + \mu^2 J^2 (g^0) \le  \delta_0^2 R^2  $$
of Lemma \ref{global.lemma} holds, with
$\delta_0 \le {1 \over 20} $ given as there.
Suppose moreover that
\begin{equation}\label{Ibound.equation}
 \lambda^2 I^2 (f^0) \le \delta_1^2 R_I^2 ,
 \end{equation}
 \begin{equation}\label{further bound.equation}
 2 \mu^2 \Gamma (2 \delta_0 R/\mu )^{2(q-1) \over q } \le \delta_1^2 R_I^2 , \
 2 \mu^2 \Gamma^q / R_I^{2-q}  \le \delta_1^2 
 \end{equation}
 and
$ {  \delta_1^2 } \le {(1- \gamma^2) \over 36} $. 
 Then
 on ${\cal T}(R) \cap {\cal T}_I(R_I)$ it holds that
 $\tau_I (\hat f - f^0 ) \le R_I $.
  \end{lemma}
 
 {\bf Proof.} We use the Basic Inequality
 $$ \| Y - \hat f - \hat g \|_n^2 + \lambda^2 I^2 (\hat f) + \mu^2 J^q (\hat g)  \le
 \| Y - f^0 - ( \hat g + \hat f_{\rm P} - f_{\rm P}^0) \|_n^2 $$ $$+ \lambda^2 I^2 (f^0) + \mu^2 J^q ( \hat g + \hat f_{\rm P} - f_{\rm P}^0 )  . $$
 This gives that
 $$ \| \hat f_{\rm A} - f_{\rm A}^0 \|^2 + \lambda^2 I^2 (\hat f ) + \mu^2 J^q ( \hat g ) $$ $$  \le
 2 \epsilon^T ( \hat f_{\rm A} - f_{\rm A}^0)/n - 2 ( \hat f_{\rm A} - f_{\rm A}^0)^T ( \hat g - g^0 + \hat f_{\rm P} - f_{\rm P}^0 ) /n$$ $$ + \| \hat f_{\rm A} - f_{\rm A}^0 \|^2 - \| \hat f_{\rm A } - f_{\rm A}^0 \|_n^2 +  \lambda^2 I^2 (f^0) + \mu^2 J^q ( \hat g + \hat f_{\rm P} - f_{\rm P}^0 ) .  $$
 By convexity the inequality also holds if we replace $\hat f$ by 
 $\tilde f := t \hat f + (1-t) \hat f^0 $
 with
 $$t:= { R_I \over R_I + \tau_I ( \hat f - f^0 ) } . $$
 Before exploiting this, we derive a bound for $J^q ( \hat g + \tilde f_{\rm P} - f_{\rm P}^0)$. 
 We use that for positive $a$ and $b$,
 $$(a+b)^q - a^q \le 2 (a+b)^{q-1} b \le 2(a^{q-1} + b^{q-1} ) b= 2a^{q-1} b + 2b^{q} $$
 Hence
 $$J^q ( \hat g + \tilde f_{\rm P} - f_{\rm P}^0) - J^q (\hat g) \le   2 J^{q-1} (\hat g ) J( \tilde f_{\rm P}  - f_{\rm P}^0 ) + 2 J^{q} (\tilde f_{\rm P} - f_{\rm P}^0 ) $$
 $$ \le  2 \Gamma J^{q-1} (\hat g ) \| \tilde f - f^0 \|  + 2 \Gamma^q \| \tilde f - f^0  \|^q $$
 where in the last step we used Condition \ref{Gamma.condition}. 
 On ${\cal T}(R)$ we have $J ( \hat g ) \le (2 \delta_0 R / \mu)^{2/q}  $ by Lemma \ref{global.lemma}. We also have 
  $\| \tilde f - f^0 \| \le R_I$. 
 Hence
 $$J^q ( \hat g + \tilde f_{\rm P} - f_{\rm P}^0) - J^q (\hat g) \le 
 2 \Gamma (2 \delta_0 R / \mu)^{2(q-1) \over q}  R_I + 2 \Gamma^q R_I^q . $$
 But then by condition (\ref{further bound.equation})
 $$\mu^2 J^q ( \hat g + \tilde f_{\rm P} - f_{\rm P}^0) \le 2 \delta_1^2 R_I^2 . $$
We insert this result in the Basic Inequality with $\hat f$ replaced by $\tilde f$:
 $$\| \tilde  f_{\rm A} - f_{\rm A}^0 \|^2 + \lambda^2 I^2 ( \tilde f) \le 6 \delta_1^2 R_I^2 +  \lambda^2 I^2 (f_0) . $$
 Invoking (\ref{Ibound.equation}) we get 
 $$\| \tilde  f_{\rm A} - f_{\rm A}^0 \|^2 + \lambda^2 I^2 ( \tilde f- f_0) \le 6  \delta_1^2 R_I^2 +
 3 \lambda^2 I^2 (f^0) \le 9 \delta_1^2 R_I^2 . $$
  Since by Lemma \ref{anti.lemma} $\| \tilde f - f^0 \|^2 \ge \| \tilde f_A - f_A^0 \|^2 / (1- \gamma^2 ) $, 
 this implies
 $$\tau_I^2 ( \tilde f - f^0 ) \le 9 \delta_1^2 R_I^2 / (1- \gamma^2 ) \le R_I^2 / 4 $$
 using $\delta_1^2 \le { (1- \gamma^2 ) \over 36 }$.
 This implies $\tau_I( \hat f - f^0 ) \le R_I $. 
 
 \hfill $\sqcup \mkern -12mu \sqcap$

\subsection{Results from empirical process theory}\label{empirical.section}

We use Theorem 2.1 in \cite{vdG2013causal} which is a consequence of
results in \cite{guedon2007subspaces} and combine this with Theorem 3.1. in \cite{vdG2013causal}.
We recall definition (\ref{Jinfty.definition}) of the entropy integral ${\cal J}_n $.
Throughout, $C_0 $ and $C_1 $ are universal constants.

\begin{theorem} \label{product2.theorem} 
Fix some $R_1$, $M_1$, $R_2$ and $M_2$ and let
$$ K_1 := \sup_{f \in {\cal F} (R_1 , M_1) }\| f \|_{\infty} , \ K_2 := \sup_{g \in {\cal G} (R_2 , M_2 ) } \| g\|_{\infty} . $$
Define for all $t$ and $n$ 
$$
B_{1,1} (t,n)  :=   { R_1  {\cal J }_n (K_1, {\cal F}(R_1 , M_1)  ) + R_1 K_1 \sqrt t   \over  \sqrt n  }  + {   {\cal J}_n^2 (K_1, {\cal F} (R_1 , M_1) ) + K_1^2 t  \over  n}    ,
 $$
 $$
B_{2,2} (t,n) :=    { R_2  {\cal J}_n (K_2,  {\cal G} (R_2 , M_2) )  + R_2 K_2 \sqrt t  \over  \sqrt n  }  + {   {\cal J}_n ^2 (K_2 , {\cal G}
 (R_2 , M_2) ) + K_2^2 t  \over  n} 
  $$
   and
 $$  B_{1,2} (t,n) := 
 { R_1 {\cal J}_n( K_2 , {\cal G}  (R_2 , M_2) )  +  R_2 {\cal J}_n  ( R_1 K_2  / R_2, {\cal F} (R_1 , M_1) )  \over  \sqrt n} 
 $$
  $$ +{ R_1 K_2 \sqrt t \over  \sqrt n} 
 +  { K_1 K_2  t \over n } . $$
 We have for all $t>0$ with probability at least $1 - C_0 \exp [-t]$
$$  \sup_{f \in {\cal F} (R_1 , M_1) } \biggl | \| f \|_n^2 - \| f \|^2 \biggr |  \le  C_1 B_{1,1} (t,n), \ 
\sup_{g \in {\cal G} (R_2 , M_2) } \biggl | \| g \|_n^2 - \| g \|^2 \biggr |  \le  C_1 B_{2,2} (t,n) . $$
Moreover, for  $R_1 K_2 \le R_2 K_1 $ and
 all values of $t $ and $n$ satisfying
  $$ C_1 B_{1,1} (t,n)  \le {R_1^2  }  , \ 
  C_1 B_{2,2} (t,n)  \le {R_2^2  }  $$
  we have with probability at least $1- C_0 \exp[-t] $
$$ \sup_{f \in {\cal F} (R_1 , M_1)  , \ g \in {\cal G}(R_2 , M_2) } \biggl | (P_n- P) fg\biggr |   \le 
C_1 B_{1,2} ( t,n)   . $$  
\end{theorem}
 
The next result follows from standard arguments using Dudley's
results (\cite{dudley1967sizes}), see e.g.\ \cite{vanderVaart:96}. 
 
\begin{theorem} \label{subGaussianproduct.theorem} 
Assume Condition \ref{subGaussian.condition} on the noise. Fix some $R_1$, $M_1$, $R_2$ and $M_2$ and let
$$ K_1 := \sup_{f \in {\cal F} (R_1 , M_1) }\| f \|_{\infty} , \ K_2 := \sup_{g \in {\cal G} (R_2 , M_2 ) } \| g\|_{\infty} . $$
 Consider values of $t$ and $n$ such that
  $$ C_1 B_{1,1} (t,n)  \le {R_1^2  }  , \ 
  C_1 B_{2,2} (t,n)  \le {R_2^2  }  $$
  with $B_{1,1} (t,n)$ and $B_{2,2} (t,n) $ given in Theorem \ref{product2.theorem}.  For these values, with probability at least $1- C_0 \exp[-t] $, one has
 $$  \sup_{f \in {\cal F} (R_1, M_1) } |  \epsilon^T f | /n\le C_1 B_{1, \epsilon} (t,n), \ 
 \sup_{g \in {\cal G} (R_2, M_2) }  |   \epsilon^T g  |/n \le C_1 B_{2, \epsilon} (t,n) , $$
 where
$$B_{1, \epsilon} (t,n) := 
 {   K_{\epsilon}  {\cal J} (R_1  , {\cal F} (R_1,M_1) )  + K_{\epsilon} R_1  \sqrt t \over  \sqrt n  }   
 $$
  and 
    $$B_{2, \epsilon} (t,n) := 
 {   K_{\epsilon}  {\cal J} (R_2 , {\cal G} (R_2,M_2) )  + K_{\epsilon} R_2  \sqrt t \over  \sqrt n  }   .
  $$
  
\end{theorem}

\begin{corollary} \label{technical.corollary} Suppose Conditions \ref{entropy.condition}, \ref{supnorm.condition} and 
\ref{subGaussian.condition}. Assume $R_1 \le M_1 / B$ and $R_2 \le M_2 / B$
 where
the constant $B$ is from Condition \ref{supnorm.condition}.
Let $B_{1,1}$, $B_{2,2}$, $B_{1,2}$ be defined as in
Theorem \ref{product2.theorem} and $B_{1 , \epsilon} $ and $B_{2, \epsilon}$ be defined
as in Theorem \ref{subGaussianproduct.theorem}.
Then 
$$ B_{1,1} (t,n)  \le    { R_1  M_1(A_I+ \sqrt t )  \over  \sqrt n  }  + {   M_1^2 (A_I^2+ t ) \over  n} , $$
$$  B_{2,2} (t,n) \le { R_2 M_2 (A_J+ \sqrt t)  \over \sqrt n} + { M_2^2  (A_J^2+ t ) \over n } . $$
$$B_{1,2} (t,n) \le 
{ A_J R_1 M_2   +  A_I R_2^{\alpha}  R_1^{1- \alpha} M_1^{\alpha} M_2^{1- \alpha}  \over  \sqrt n} 
 +{ R_1 M_2 \sqrt t \over  \sqrt n} 
 +  { M_1 M_2  t \over n }  $$
 and
$$  B_{1, \epsilon} (t,n) \le  
 {   A_I K_{\epsilon}  M_1^{\alpha} R_1^{1-\alpha}  + K_{\epsilon} R_1  \sqrt t \over  \sqrt n  }  , \ 
 B_{2, \epsilon} (t,n) \le
 {   A_J K_{\epsilon}  M_2^{\beta} R_2^{1-\beta}  + K_{\epsilon} R_2  \sqrt t \over  \sqrt n  } . $$
  The constants $A_I$ and $A_J$ are from Condition \ref{entropy.condition} and
 the constant $K_{\epsilon}$ from Condition \ref{subGaussian.condition}.
\end{corollary}

\begin{theorem} \label{technical.theorem}
Assume Conditions  \ref{entropy.condition}, \ref{supnorm.condition} and 
\ref{subGaussian.condition}.
Let $\lambda \le  R_I \le \mu \le R \le 1$ be constants
and $L_I := R_I / \lambda $ and $L_J := (R/\mu)^{2/q} $.\\
{\bf Case 1.} Assume $\lambda^2 \le 1/ B^2 $ and $\mu^2 \le R^{2-q } / B^{q}  $.
Suppose that for some
$L \ge 4 C_1$, 
\begin{equation}\label{lambdamu.equation}
\sqrt n \lambda^{1+ \alpha}  \ge   L   A_I, \ 
\sqrt n \mu^{1+ \beta}  \ge L  A_J , 
\end{equation}
\begin{equation}\label{R.equation}
 R \ge L L_J A_J / \sqrt n , \ R \ge K_{\epsilon} \lambda , \ R \ge L_J  \lambda , \ R \ge K_{\epsilon }^{q  \over q-
 (2-q) \beta }  \mu   
 \end{equation}
and
\begin{equation}\label{lambda.equation}
 \lambda^{\alpha}  \le 1/L .
 \end{equation}
Then with probability at least $1-3 C_0 \exp[-n \lambda^2 / L^2 ]$ it holds that
$$ \sup_{f \in {\cal F} (R, R/\lambda) }\biggl |  \| f \|_n^2 - \| f \|^2 \biggr | \le {4 R^2 \over L} , \ 
 \sup_{g \in {\cal G} (R, L_J) } \biggl | \| g \|_n^2 - \| g \|^2 \biggr | \le {4 R^2 \over L } , $$
$$ \sup_{f \in {\cal F} (R , R/\lambda) , \ g \in {\cal G} (R, L_J) } \biggl | 
( P_n - P) fg \biggr | \le {4 R^2 \over L} $$
and 
$$ \sup_{f \in {\cal F} (R, R/\lambda) }  |   \epsilon^T f   |/n \le {2 R^2 \over L} , \ 
\sup_{g \in {\cal G} (R,L_J ) } |  \epsilon^T g   | /n\le {2 R^2 \over L } . $$
{\bf Case 2.} Assume moreover that
\begin{equation} \label{RI.equation}
R_I \ge L L_J A_J / \sqrt n, \  R_I \ge K_{\epsilon} \lambda
\end{equation}
Then with probability at least $1- 3 C_0 \exp[-n \lambda^2 / L^2 ] $,
$$\sup_{f \in {\cal F} (R_I, L_I ) } \biggl |  \| f \|_n^2 - \| f \|^2 \biggr | \le {4 R_I^2 \over L} , \ 
\sup_{f \in {\cal F} (R_I , L_I  ) , \ g \in {\cal G} (R,L_J) } \biggl | 
( P_n - P) fg \biggr | \le {4 R_I^2 \over L} $$
and
$$ \sup_{f \in {\cal F} (R_I , L_I ) }  |   \epsilon^T f   |/n \le {2 R_I^2 \over L} .$$
\end{theorem}

{\bf Proof of Theorem \ref{technical.theorem}.}

{\bf Case 1.}
We first apply Corollary \ref{technical.corollary} 
with $R_1=R_2=R$ and $M_1= R/\lambda$, $M_2=L_J := (R/\mu)^{2/q}$.
The condition $\lambda \le 1/B $ ensures $R_1 \le M_1/ B$ and the condition
$\mu \le R^{2-q \over q} / B$ ensures that $R_2 \le M_2 / B $. 
We let $B_{1,1}$, $B_{2,2}$, $B_{1,2}$ be defined as in
Theorem \ref{product2.theorem} and $B_{1 , \epsilon} $ and $B_{2, \epsilon}$ be defined
as in Theorem \ref{subGaussianproduct.theorem} and insert the value $t= n \lambda^2 / L^2 $.

{\bf Case 1a for $\| f\|_n^2$.}
$$B_{1,1} (t,n)  
 \le  { R^2   (A_I + \sqrt t )  \over  \sqrt n \lambda   }  + {  R^2 (A_I^2+  t ) \over  n \lambda^2}   
 =  \biggl ( {    A_I + \sqrt t  \over  \sqrt n \lambda   }  + {   A_I^2+  t  \over  n \lambda^2}   \biggr ) R^2 . $$
 Now use that by (\ref{lambdamu.equation}) $\sqrt n \lambda  \ge L  A_I $ and $t = n \lambda^2/L^2  $ to get
 $$B_{1,1} (t,n) \le\biggl (  {2 \over L} + {2 \over L^2}\biggr )  \le {4 R^2  \over L} $$
 
 {\bf Case 1b For $\| g\|_n^2$.}
 $$B_{2,2} (t,n) 
\le 
  { R L_J (A_J   +  \sqrt t)   \over  \sqrt n  }  + {  L_J^2  (A_J^2   + t )  \over  n} 
    = 
   \biggl (  {  L_J (A_J  +  \sqrt t )   \over   \sqrt n R  }  + {  L_J^2 ( A_J^2  + t )  \over n R^2 } 
   \biggr )R^2
    $$
      $$ \le \biggl ( { 1 \over L} +{  \sqrt {t } \over \sqrt n R } + {1 \over L^2} + { t \over n R^2 } \biggr ) R^2 ,
  $$
  where we used that $R \ge L L_J A_J / \sqrt n $ by (\ref{R.equation}). Insert $R \ge \lambda$ and
  $t = n \lambda^2 / L^2 $ to get
  $$B_{2,2} (t,n) \le \biggl ( {2 \over L} + { 2 \over L^2} \biggr ) R^2 \le { 4 R^2 \over L }  . $$

{\bf Case 1c for $f^Tg/n$.}
We already know by Cases 1a and 1b that 
$C_1 B_{1,1} (t,n) \le R^2 $ and $C_1 B_{2,2} (t,n) \le R^2 $ with probability
at least $1- C_0 \exp[- n \lambda^2 / L^2 ] $.  Moreover
   $$
  B_{1,2} (t,n) 
  \le  
   { R  L_J A_J  +  R^{1+ \alpha}  L_J^{1- \alpha} A_I  / \lambda^{\alpha}   + R  L_J \sqrt t \over  \sqrt n} 
    +  { R  L_J t \over n  \lambda }  $$
    $$ =  \biggl ( { L_J A_J  \over  \sqrt n R}    + {    L_J^{1-\alpha} A_I  \lambda  \over  \sqrt n R^{1 - \alpha}  \lambda^{1+\alpha} }   + {   L_J \sqrt t \over   \sqrt n R } 
    +  {  L_J  t  \over n R  \lambda   }  \biggr ) { R^2  }   . $$
    Use $R \ge \lambda L_J $, $R \ge L L_J A_J /\sqrt n $  from
    (\ref{R.equation}) and $ \sqrt n \lambda^{1+\alpha} \ge L A_I
    $ from (\ref{lambdamu.equation}) to find that 
  $$ B_{1,2} (t,n)  \le  \biggl ( { 1 \over L }    +  \lambda^{\alpha}   + {    \sqrt t \over \sqrt n  \lambda } 
    +  {   t  \over  n  \lambda^2   }  \biggr ) { R^2 }   . $$
    Apply now that by (\ref{lambda.equation}) $\lambda^{\alpha} \le 1/L$ and $t = n \lambda^2 /L^2$  to get
    $$B_{1,2} (t,n) \le \biggl ( {3 \over L} + {1 \over L^2} \biggr ) \le  { 4 R^2 \over L } . $$

{\bf Case 1d for $\epsilon^T f/n $.}
We already know by Cases 1a and 1b that 
$C_1 B_{1,1} (t,n) \le R^2 $ and $C_1 B_{2,2} (t,n) \le R^2 $ with probability
at least $1- C_0 \exp[- n \lambda^2 / L^2 ] $. Moreover
$$B_{1, \epsilon} (t,n) \le { K_{\epsilon}A_I  \lambda R  \over \sqrt n \lambda^{1+ \alpha} } +
{ K_{\epsilon} R \sqrt t \over \sqrt n } 
=  \biggl ( { K_{\epsilon}A_I  \lambda   \over \sqrt n \lambda^{1+ \alpha} R } +
{ K_{\epsilon} \sqrt t \over \sqrt n  R} \biggr ) R^2  .$$
Invoke $\sqrt n \lambda^{1+ \alpha} \ge L A_I  $ from (\ref{lambdamu.equation}) and $R \ge K_{\epsilon} \lambda  $ 
from (\ref{R.equation}) to
obtain
$$ B_{1, \epsilon} (t,n) \le \biggl ( {1 \over L} + {\sqrt t \over \sqrt n \lambda } \biggr ) R^2 . $$
With $t= n \lambda^2 /L^2 $ this gives
$$ B_{1, \epsilon} (t,n) \le {2 R^2 \over L } . $$

{\bf Case 1e for $\epsilon^T g/n  $.}
We gave
$$ B_{2, \epsilon} (t,n)  \le  \biggl ( { K_{\epsilon}   A_J  L_J^{\beta} \over 
\sqrt n R^{1+ \beta} } + {K_{\epsilon} \sqrt t  \over \sqrt n R } \biggr ) R^2  .$$
Use $\sqrt n \mu^{1 + \beta}  \ge L  A_J $ from (\ref{lambdamu.equation})  to find
$$ B_{2, \epsilon} (t,n)\le \biggl ( { L_J^{\beta} ( \mu / R )^{1+\beta} K_{\epsilon}  \over L}  +
{K_{\epsilon} \sqrt t  \over \sqrt n R }
\biggr ) R^2 .$$ 
Next, we see that $ L_J^{\beta} ( \mu / R )^{1+\beta} \le 1/ K_{\epsilon} $ since $R \ge K_{\epsilon}^{q \over q-(2-q)\beta } \mu$ by
(\ref{R.equation}). Moreover, also by (\ref{R.equation}) $R \ge K_{\epsilon} \lambda $.
So with $t = n \lambda^2 / L^2 $
$$ B_{2, \epsilon} (t,n) \le \biggl ( {1 \over L} + {\sqrt t \over \sqrt n \lambda } \biggr ) R^2 =
{2 R^2 \over L } . $$

{\bf Case 2.} Take $R_1=R_I$, $R_2=R$ and $M_1= L_I $, $M_2 = L_J $. 
Then again 
$R_1 \le M_1 / B$ and $R_2 \le M_2 / B$. Also With these new values,
we let $B_{1,1}$, $B_{2,2}$, $B_{1,2}$ be defined as in
Theorem \ref{product2.theorem} and $B_{1 , \epsilon} $ and $B_{2, \epsilon}$ be defined
as in Theorem \ref{subGaussianproduct.theorem} and insert the value $t= n \lambda^2 / L^2 $.

{\bf Case 2a for $\| f\|_n^2$.}
$$B_{1,1} (t,n) \le { R_I^2   ( A_I  + \sqrt t ) \over \sqrt n \lambda}+
R_I^2 { A_I^2 + t \over n \lambda^2 } = \biggl ( { A_I +  \sqrt t \over \sqrt n R_I } +
{ A_I + t \over nR_I^2 } \biggr ) R_I ^2  \le {4 R_I^2 / L } $$
since $\sqrt n R_I \ge \sqrt n \lambda \ge \sqrt n \lambda^{1+ \alpha} \ge L A_I $ by
(\ref{lambdamu.equation}) and
$t = n \lambda^2 / L^2 $. 

{\bf Case 2b for $f^Tg/n$.}
By similar arguments as in Case 1a (see also Case 2a) and 1b that 
$C_1 B_{1,1} (t,n) \le R_I^2 $ and $C_2 B_{2,2} (t,n) \le R^2 $ with probability
at least $1- C_0 \exp[- n \lambda^2 / L^2 ] $.
Moreover
$$B_{1,2} (t,n) \le { R_I L_J A_J \over \sqrt  n } + {R^{\alpha} R_I \lambda L_J^{1- \alpha}  A_I \over \sqrt n 
\lambda^{1+ \alpha} } +
{ R_I L_J \sqrt t \over \sqrt n} + {t R_I L_J \over n \lambda } $$
$$ = \biggl (  {  L_J A_J \over \sqrt  n R_I  } + {R^{\alpha} \lambda L_J^{1- \alpha}  A_I \over \sqrt n 
\lambda^{1+ \alpha} R_I  } +
{  L_J \sqrt t \over \sqrt n R_I } + {t L_J \over n \lambda R_I  }   \biggr ) R_I^2 .
 $$
 Use that
 $R_I \ge L L_J A_J / \sqrt n$ (see (\ref{RI.equation})), $ \sqrt n \lambda^{1+ \alpha} \ge L A_I$  
 (see (\ref{lambdamu.equation})) and
 $ R \ge \lambda L_J$ (see (\ref{R.equation})).  We then get
 $$B_{1,2} (t,n) \le \biggl ( {1 \over L} + \lambda^{\alpha} + {\sqrt t \over
 \sqrt n \lambda } + {t \over n \lambda^2 } \biggr ) R_I^2 . $$
 With $t = n \lambda^2 / L^2$  and $\lambda^{\alpha } \le 1/ L $ (see (\ref{lambda.equation})) this gives again
 $$B_{1,2} (t,n)  \le { 4 R_I^2 \over L } . $$

{\bf Case 2c for $\epsilon^T f/n $.}
By Case 2a, it holds that
$C_1 B_{1,1} (t,n) \le R_I^2 $ with probability
at least $1- C_0 \exp[- n \lambda^2 / L^2 ] $. Moreover
$$B_{1, \epsilon}  (t,n) \le  \biggl (  { K_{\epsilon} A_I  \lambda \over
\sqrt n \lambda^{1+ \alpha} R_I } +  { K_{\epsilon}  \sqrt t \over \sqrt n R_I } \biggr ) R_I^2  . $$
From (\ref{lambdamu.equation}) we know $\sqrt n \lambda^{1+ \alpha} \ge L A_I$ and from
(\ref{RI.equation}) 
$R_I \ge K_{\epsilon } \lambda$. With $t = n \lambda^2 / L^2$ we find
$$B_{1, \epsilon}  (t,n) \le { 2 R_I^2 \over L }  . $$
The result now follows from the same arguments as in Case 2 of Theorem \ref{technical.theorem}.

\hfill $\sqcup \mkern -12mu \sqcap$

\begin{remark} \label{simplify.remark} If we assume condition (\ref{lambdamu.equation}), then condition (\ref{R.equation}) is met for 
$$ K_{\epsilon} { \lambda \over \mu } \le K_{\epsilon}^{q \over q- (2-q) \beta } \le
{R \over \mu } \le {\rm min}\biggl  \{ \biggl ({ \sqrt n \over L A_J  } \biggr  )^{ q \beta \over (1+ \beta) (2- q )} ,
\biggl ( {\mu \over \lambda} \biggr  )^{q \over 2- q } \biggr \} . $$
Under general conditions, the left hand side tends to zero  and
the right hand side tends to infinity as $n \rightarrow \infty$. 

\end{remark}

\subsection{Application to ${\cal T } (R) $}\label{TR.section}

Recall the definition (\ref{TR.equation}) of the set ${\cal T} (R)$.

\begin{lemma} \label{TR.lemma} Let $\lambda \le \mu \le R \le 1$.
Assume Conditions  \ref{entropy.condition}, \ref{supnorm.condition}, 
\ref{subGaussian.condition} and \ref{incoherence.condition}.
Assume that $\lambda^2 \le (1- \gamma)/B^2 $ and $\mu^2 \le (1- \gamma)^q R^{2-q} / B^q $.
Let
$$ L \ge \max \biggl \{ 4 C_1 (1- \gamma)^{1/2}, 16 / ((1- \gamma)^{1/2} \delta_0^2 )\biggr  \}  $$
and
$$ \sqrt n \lambda^{1+ \alpha}  \ge   L   A_I, \  \sqrt n \mu^{1+ \beta}  \ge L 
A_J , $$
$$\lambda^{\alpha} \le (1- \gamma)^{1+\alpha \over 2} /L $$
and
$$ K_{\epsilon} { \lambda \over \mu } \le K_{\epsilon}^{q \over q- (2-q) \beta } \le
{R \over \mu } \le {\rm min}\biggl  \{ \biggl ({ \sqrt n (1- \gamma)^{1/2}\over L A_J  } \biggr  )^{ q \beta \over (1+ \beta) (2- q )} ,
\biggl ( {\mu \over \lambda} \biggr  )^{q \over 2- q } \biggr \} . $$
Then $$\PP ({ \cal T} (R)  ) \ge1- \exp[ - n  \lambda^2 / L^2 ] . $$
\end{lemma}

{\bf Proof.} Recall the definition of ${\cal M}(R)$ given in (\ref{MR.equation}) with $\tau_R( \cdot , \cdot) $ given in
(\ref{tau.equation}). 
Define $\tilde \lambda^2 := \lambda^2 / (1- \gamma)$,
$\tilde \mu^2 := \mu^2 / (1- \gamma)$ and  $\tilde R^2 := R^2 / (1- \gamma) $. 
By Lemma \ref{incoherence.lemma}
$${\cal M}(R) \subset \biggl \{ (f,g):\ \| f \| \le \tilde R , \ \| g \| \le \tilde R ,\  I(f) \le R/ \lambda , \  
J(g) \le (R / \mu )^{2 \over q}  \biggr \} $$
$$ = \biggl \{ (f,g): \ f \in {\cal F} ( \tilde R , \tilde R / \tilde \lambda) , \ g \in {\cal G} ( \tilde R , 
(\tilde R / \tilde \mu )^{2 \over q}\biggr  \} . $$
We apply Case 1 of Theorem \ref{technical.theorem} with $(\lambda , \mu , R)$ replaced by $
(\tilde \lambda , \tilde \mu , \tilde R)$. We also replace
$L$ by $\tilde L ^2:= L^2 /(1- \gamma)$.
Then 
$$  \sqrt n \tilde \lambda^{1+ \alpha}= \sqrt n \lambda^{1+ \alpha} /(1-\gamma)^{1+ \alpha \over 2} $$
$$ \ge L A_I / (1-\gamma)^{1+ \alpha \over 2}= \tilde L A_I  /(1-\gamma)^{\alpha \over 2} \ge
\tilde L A_I . $$
Similarly
$$\sqrt n \tilde \mu^{1+\beta} \ge \tilde L A_J . $$
The condition $\lambda^{\alpha} \le (1- \gamma)^{1+\alpha \over 2} /L $ gives
$$ \tilde \lambda^{\alpha} = { \lambda^{\alpha} \over (1- \gamma)^{\alpha / 2} } \le
{ (1- \gamma)^{1+ \alpha \over 2} \over L  (1- \gamma)^{\alpha \over 2}} = {1 \over  \tilde L }  .$$
Furthermore
$$ K_{\epsilon} { \tilde \lambda \over \tilde  \mu } \le K_{\epsilon}^{q \over q- (2-q) \beta } \le
{\tilde R \over \tilde \mu } \le {\rm min}\biggl  \{ \biggl ({ \sqrt n \over \tilde L A_J  } \biggr  )^{ q \beta \over (1+ \beta) (2- q )} ,
\biggl ( {\tilde \mu \over\tilde  \lambda} \biggr  )^{q \over 2- q } \biggr \} . $$
By Remark \ref{simplify.remark} we conclude that the conditions for Case 1 of
Theorem \ref{technical.theorem} are met.
Clearly, for any $f$ and $g$
$$\biggl | \| f+ g \|_n^2 - \| f+ g \|^2 \biggr | \le \biggl |  \| f \|_n^2 - \| f \|_n^2 \biggr | +\biggl |  \| g \|_n^2 - \| g \|^2 
\biggr | +
\biggl | 2 (P_n - P )fg \biggr | . $$
By Case 1 of Theorem \ref{technical.theorem}, for 
$\tilde L = L/(1- \gamma)^{1/2} \ge 4 C_1$ and 
for $16 \tilde R^2 / \tilde L \le \delta_0^2 R^2$
$$\PP ({ \cal T} (R)  ) \ge 1- \exp[ - n \tilde \lambda^2 / \tilde L^2 ] . $$
The proof if finished by noting that $\tilde R^2 / \tilde L = R^2 / ( L(1- \gamma)^{1/2} )$
and $\tilde \lambda^2 / \tilde L^2 = \lambda^2 / L^2 $.

 \hfill $\sqcup \mkern -12mu \sqcap$

 \subsection{Application to ${\cal T}_I (R_I) $}\label{TRII.section}
 
 Recall the definition (\ref{TRI.equation}) of the set ${\cal T}_I (R_I, R)$.

 \begin{lemma} \label{TRI.lemma} 
 Assume Conditions  \ref{entropy.condition}, \ref{supnorm.condition}, 
\ref{subGaussian.condition}, \ref{incoherence.condition} and \ref{Gamma.condition}.
 Let $\lambda\le R_I  \le \mu \le R \le 1$. 
 Assume that $\lambda^2 \le (1- \gamma)/(2B)^2 $ and $\mu^2 \le (1- \gamma)^q R^{2-q} / (2B)^q $.
Let
$$ L \ge \max \biggl \{ 2 C_1 (1- \gamma)^{1/2}, 32 / ((1- \gamma)^{1/2} \delta_0^2 ),
32 / ( \delta_1^2 ) \biggr  \}  .$$
Take
$$ \sqrt n \lambda^{1+ \alpha}  \ge   L   A_I, \  \sqrt n \mu^{1+ \beta}  \ge L 
A_J , $$
$$\lambda^{\alpha} \le (2(1- \gamma))^{1+\alpha \over 2} /L $$
and
$$ K_{\epsilon} { \lambda \over \mu } \le K_{\epsilon}^{q \over q- (2-q) \beta } \le
{R \over \mu } \le {\rm min}\biggl  \{ \biggl ({ \sqrt n (1- \gamma)^{1/2}\over 2 L A_J  } \biggr  )^{ q \beta \over (1+ \beta) (2- q )} ,
\biggl ( {\mu \over \lambda} \biggr  )^{q \over 2- q } \biggr \} . $$
Also take
$$ R_I \ge L ( R / \mu)^{2 \over q}  A_J / \sqrt n, \  R_I \ge K_{\epsilon} \lambda $$
$$ \Gamma R_I \le  (2 R/\mu)^{2 \over q} $$
Then
$$\PP ( {\cal T}_I (R_I, R))  \ge 1- 3C_0 \exp[- n \lambda^2 / L^2 ] .$$

  \end{lemma} 
 
 {\bf Proof.} By Lemma \ref{entropy.lemma}
 $${\cal J}_{\infty} (z , \{f_{\rm A} : \ f \in {\cal F} (R,M) \} ) \le 2 {\cal J}_{\infty} (z, {\cal F} (R,M) ), \ z > 0  $$
and for $R \le M / B$
$$\sup_{f \in {\cal F}(R,M) } \| f_{\rm A} \|_{\infty} \le 2  M . $$
We can therefore apply similar arguments as for Case 2 of Theorem \ref{technical.theorem}.
We know that for $f \in {\cal F} ( R_I)$, $\| f_{\rm A} \| \le \|f \| \le R_I $. So
$${\cal J}_{\infty} (z , \{f_{\rm A} : \ f \in {\cal F} (R_I) \} ) \le 2 {\cal J}_{\infty} (z, {\cal F} (R_I ,R_I / \lambda ) ), \ z > 0  $$
and
$$\sup_{f \in {\cal F}(R_I) } \| f_{\rm A} \|_{\infty} \le 2  R_I / \lambda  . $$
Moreover, for $f \in {\cal F} (R_I)$ and $g \in {\cal G} (R)$ we have
$$J( g+ f_{\rm P} ) \le J(g) + J( f_{\rm P}) \le (R/\mu)^{2/q} + \Gamma \| f \| \le
(R/\mu)^{2/q} + \Gamma R_I \le  (2 R/\mu)^{2/q} , $$
and 
$$ \| g + f_{\rm P} \| \le \| g \| + \| f_{\rm P} \| \le R / (1-\gamma)^{1/2} + R_I \le 2 R/(1-\gamma)^{1/2}  . $$
It follows that
$$ \{ g+ f_{\rm P} : \ f \in {\cal F} (R_I) , \ g \in {\cal G} (R) \} \subset
{\cal G} ( 2R/ (1- \gamma)^{1/2} , (2 R/\mu )^{2/q } ) . $$
It is also clear that for any $f$ and $g$
$$P f_{\rm A} g =  E  f (X_1) g(Z_1) - E \biggl [  E (f(X_1 ) \vert Z )  g(Z_1) \biggr ] =0 $$
and similarly $P f_{\rm A} f_{\rm P} = 0 $. 
By an appropriate replacements of the constants in Case 2 of Theorem \ref{technical.theorem}
(as in the proof of Lemma (\ref{TR.lemma}) now using $(1- \gamma)^{1/2} / 2$ instead of
$(1- \gamma)^{1/2}$) the results follows. 
 
  \hfill $\sqcup \mkern -12mu \sqcap$

\subsection{Finishing the proof of Theorem \ref{main.theorem}} \label{finish.section}
We first note that 
since $\max\{ A_I, A_J \} \le n^{1- \delta \over 2} $
we $\lambda^{1+ \alpha} = c_1 A_I/\sqrt n \le n^{-\delta/2}$.
So for $n$ large $\lambda$ will be small. The same is true for $\mu$ and
for the ratio $\lambda / \mu$.

In view of Lemma \ref{global.lemma} we need $\lambda^2 I^2 (f^0) +
\mu^2 J^q (g^0) \le \delta_0^2 R^2 $. We take 
$$R^2 =\max \biggl \{  \mu^2 J^q (g^0 ) / (4\delta_0^2 ) ,  K_{\epsilon}^{2 q \over 2- (2-q) \beta } \biggr \} $$
and $n $ sufficiently large such that
$$\lambda^2 I^2 (f^0) \le \mu^2 J^q (g^0) .$$
Take
$$ L = \max \biggl \{ 2 C_1 (1- \gamma)^{1/2}, 32 / ((1- \gamma)^{1/2} \delta_0^2 ),
32 / ( \delta_1^2 ) \biggr  \}  .$$
Since
$$\max \biggr \{  { J^{q/2} (g^0) \over 2 \delta_0}  ,  K_{\epsilon}^{ q \over 2- (2-q) \beta } \biggr \} \le 
 {\rm min}\biggl  \{ \biggl ({ \sqrt n (1- \gamma)^{1/2}\over 2 L A_J  } \biggr  )^{ q \beta \over (1+ \beta) (2- q )} ,
\biggl ( {\mu \over \lambda} \biggr  )^{q \over 2- q } \biggr \}  $$
for $n$ sufficiently large as 
 $A_J \le n^{1  - \delta \over 2  }$ we know 
 from Remark \ref{simplify.remark}
that the conditions for Lemma \ref{TR.lemma}
are met for $n$ sufficiently large.
By Lemma \ref{tighter.lemma} we also need $\lambda^2 I^2 (f_0) \le R_I^2 / \delta_1^2 $.
For  $R_I / \lambda = \max \{ I(f^0) / \delta_1 ,  K_{\epsilon} \} $.
$$ R_I \ge L ( R / \mu)^{2 \over q}  A_J / \sqrt n$$
for $n$ sufficiently large so we can also apply Lemma \ref{TRI.lemma}.

\bibliographystyle{plainnat}
\bibliography{reference}

\end{document}